\documentclass[11pt]{article}
\usepackage[latin1]{inputenc}
\usepackage{amsmath}
\usepackage{float}
\usepackage{mathtools}
\usepackage{comment}
\usepackage{upgreek}
\input{insbox}

\usepackage{multicol}
\usepackage{amsfonts}
\usepackage{mathrsfs}
\usepackage{bm}
\usepackage{gensymb}
\usepackage{wrapfig}
    \input{insbox}
\usepackage{enumerate}
\usepackage{enumitem}

\usepackage{amssymb}
\usepackage{titlesec}
\usepackage{marvosym}
\usepackage{array}
\usepackage{chngcntr}
\usepackage{blindtext}
\usepackage[english]{babel}
\usepackage{amsthm}
    \newtheoremstyle{mythm}
    {24pt} % Space above
    {24pt} % Space below
    {} % Body font
    {} % Indent amount
    {\bfseries} % Theorem head font
    {.} % Punctuation after theorem head
    {.5em} % Space after theorem head
    {} % Theorem head spec (can be left empty, meaning `normal')
    \theoremstyle{mythm}
    \newtheorem{thm}{Theorem}[section]
    
    \newtheorem{lemma}[thm]{Lemma}
    \newtheorem{remark}[thm]{Remark}

    \newtheorem{definition}[thm]{Definition}
    \theoremstyle{remark}
\usepackage{hyperref} %for generating hyperlinks
    \hypersetup{
    colorlinks=true,
    linkcolor=blue,
    filecolor=blue,      
    urlcolor=blue,
    citecolor=blue,
    }
\usepackage{graphicx}
\usepackage{caption}
\usepackage{color}
\usepackage{authblk}
\usepackage[titletoc,title]{appendix}
\usepackage{cite}
\usepackage[top=1in, bottom=1in, left=1in, right=1in]{geometry}

\def\rw{\ensuremath\rightarrow }

\def\v{\ensuremath\bigcup}

\def\0{\ensuremath\varnothing}

\makeatletter
\newcommand{\compemb}{\mathrel{\mathpalette\comp@emb\relax}}
\newcommand{\comp@emb}[2]{%
  \vcenter{%
    \offinterlineskip\m@th
    \ialign{$#1##$\cr\hookrightarrow\cr\noalign{\vskip1pt}\hookrightarrow\cr}%
  }%
}
\makeatother
\def\a{\ensuremath\alpha }

\def\g{\ensuremath\gamma }

\def\d{\ensuremath\delta }

\def\e{\ensuremath\varepsilon }
\def\f{\ensuremath\varphi }

\def\o{\ensuremath\omega }
\def\O{\ensuremath\Omega }
\def\u{\ensuremath\bm{u} }
\def\v{\ensuremath\bm{v} }
\def\x{\ensuremath\bm{x} }
\def\w{\ensuremath\bm{w} }

\def\C{\ensuremath\mathcal{C} }

\def\B{\mathscr{B} }

\def\L{\mathscr{L} }
\def\K{\ensuremath\mathcal{K} }

\def\H{\ensuremath\mathcal{H} }

\def\R{\ensuremath\mathbb{R} }

\def\N{\ensuremath\mathbb{N} }

\renewcommand{\div}{\text{div}}

\newlength{\jeroenlen}

\counterwithin*{equation}{section}

\titleformat{\section}{\bfseries\fontsize{15.5}{0}\selectfont}{\thesection.}{0.5em}{}
\renewcommand{\thesection}{\normalfont\arabic{section}}
\titlespacing*{\section}{0pt}{36pt}{6pt}

\titleformat{\subsection}{\bfseries\fontsize{12}{0}\selectfont}{\thesubsection.}{0.5em}{}
\renewcommand{\thesubsection}{\normalfont\arabic{section}.\arabic{subsection}}
\titlespacing*{\subsection}{0pt}{18pt}{3pt}

\numberwithin{equation}{section}

\numberwithin{thm}{section}

\numberwithin{figure}{section}

\numberwithin{table}{section}

\title{On the Propulsion of a Rigid Body in a Viscous Liquid by Time-Periodic Force with a Zero Average} 
\author[1]{Joris Edelmann} 
\author[2]{Giovanni P. Galdi}
\author[2]{Mher M. Karakouzian}
\author[1]{Thomas Richter}

\affil[1]{Faculty of Mathematics, Institute of Analysis and Numerics, Otto von Guericke University Magdeburg, 	Magdeburg, Germany}
\affil[2]{Department of Mechanical Engineering \& Materials Science, SSOE\newline University of Pittsburgh. Pittsburgh, USA}

\usepackage{booktabs}

\graphicspath{%
    {converted_graphics/}% inserted by PCTeX
    {/}% inserted by PCTeX
}
\begin{document}
\renewcommand{\labelenumi}{(\roman{enumi})}

    \maketitle

\begin{abstract}
We perform analytical and numerical analyses of the propulsion of a rigid body in a viscous fluid subjected to a periodic force with zero average over a period. This general formulation specifically addresses the significant case, where propulsion is generated by the oscillation of a mass located in an internal cavity of the body. We provide a rigorous proof of the necessary and sufficient conditions for propulsion at the second order of magnitude of the force. These conditions are implemented and confirmed by numerical tests for bodies without fore-and-aft symmetry, while they are silent for bodies with such symmetry, like round ellipsoids. Consequently, in this case, propulsion can only occur at an order higher than the second. This problem is investigated by numerically integrating the entire set of equations, and the result shows that, in fact, propulsion does occur, thus opening new avenues for further analytical studies.\medskip\par\noindent


\end{abstract}
\section{Introduction}\label{ssec:intro}

In the last two decades, there has been considerably growing interest in studying the propulsion of a rigid body $\mathscr B$ immersed in a resistive medium, such as a viscous liquid $\mathscr L$ where the propulsive mechanism is due to an internal mass performing oscillatory motions. A fairly comprehensive overview of the related problems, results, and various remarkable applications is available in the monograph \cite{BC} and the references therein.

This phenomenon is, at first sight surprising and definitely captivating for several reasons. In the first place, the net force (average
over a period) acting on $\mathscr B$ and due to the oscillations of the internal mass is zero, yet propelling. Moreover, propulsion seems to occur ---under suitable mass oscillations--- even when $\mathscr B$ and the internal cavity $\mathcal C$, where the mass is located,  possess full symmetry; for example, $\mathscr B$ is a spherical shell and $\mathcal C$ is its interior \cite{V}. This fact suggests that the oscillations to which $\mathscr B$ is subject must induce, somehow, an ``asymmetry" in the surrounding liquid flow that is thus able to produce an overall non-zero thrust. Therefore, the inertia of the liquid should play an important role. Restated in mathematical terms, this means that the phenomenon is genuinely nonlinear. In this regard, it should be noted that the mathematical analysis used to understand this type of propulsion has been performed, so far, on rather simplified models where the action of $\mathscr L$ on $\mathscr B$ is a prescribed viscous force depending on the velocity of the center of mass $G$ of $\mathscr B$. As a result, the motion of the body is somehow decoupled from that of the liquid and the motion of $G$ is reduced to the study of a system of ordinary differential equations; see \cite[Chapter 5]{BC}. However, this simplification can lead to neglecting some important properties due precisely to the {\em mutual}  interaction between the body and the liquid.

Very recently, we have started a rigorous mathematical analysis aimed at understanding the propulsion of a rigid body in a viscous liquid caused by a generic time-periodic force $\bm{f}$ acting on the body \cite{galdi2025propulsion}. The model we used is the complete one; namely, obtained by coupling  Newton's equations for the body with the  Navier-Stokes equations for the liquid. In  \cite{galdi2025propulsion}, we restricted our attention to the case where the average of $\bm{f}$ over a period ---denoted by $\overline{\bm{f}}$--- is not zero. In particular, we  furnished necessary and sufficient conditions on $\bm{f}$ and the physical properties of $\mathscr B$ ensuring propulsion, at least at the first order in the magnitude $\delta$ of $\bm{f}$. As somehow expected, these findings show some remarkable properties of propulsion that simplified models might not able to catch \cite[Section 5]{galdi2025propulsion}.\looseness=-1

Although relevant, however, the results obtained in \cite{galdi2025propulsion} are incomplete, because they cannot cover the situation where the propulsion is caused by an oscillating internal mass, since, in that case, we have $\overline{\bm{f}}={\bf 0}$; see \eqref{eq:firsAppofForce} and the comments that follow.

The goal of this paper is to provide a first rigorous contribution to the propulsion of a rigid body in a viscous fluid under the action of a periodic force with zero mean over a period. This allows, as a special case, propulsion by an oscillating internal mass. To simplify the scenario, we assume that $\mathscr B$ cannot rotate, which can be practically achieved by imposing a suitable torque on $\mathscr B$. Our study is supported and complemented by targeted numerical tests that, on the one hand, validate the analytical results and, on the other, predict entirely new features not captured by our current theory. As such, they open new avenues of investigation that will be explored in future work.\looseness=-1

The condition that ensures that a body subjected to an oscillating driving mechanism of period $T$ actually has a non-zero net motion is that its center of mass $G$ can travel any given finite distance in a finite time. Denoting by $\bm{\gamma}$ the velocity of $G$, it is then easily shown (see Section 2) that this happens if and only if
\begin{equation}\label{gamma}  
    \overline{\bm{\gamma}}:=\frac1T\int_0^T\bm{\gamma}(t)\,{\rm d}t\neq{\bf 0}\,.
\end{equation}
Therefore, our objective is to find conditions on $\bm{f}$ ensuring the validity of \eqref{gamma}. To this end, we write $\bm{f}=\d\bm{F}$, where $\d$ represents the  magnitude  of the oscillations (in appropriate norm), and begin by observing that, since $\overline{\bm{f}}={\bf 0}$, it follows that, at the order of $\delta$, it is  $\overline{\bm{\gamma}}={\bf 0}$; see Remark \ref{ave}. Consequently, propulsion can only occur at $O(\delta^2)$ (or higher), meaning that it is a genuine {\em nonlinear} phenomenon, as expected. Our objective is then to establish necessary and sufficient conditions ensuring \eqref{gamma} at the order of $\delta^2$. More precisely, we prove that in a rather general class of weak solutions (see Theorem \ref{thm:ExistceWkSol_NonlnrProb}), the following relation holds:
\begin{equation}\label{th}
    \overline{\bm{\gamma}}=\delta^2\,\mathbb A\cdot\bm{G}+\bm{\mathsf{R}}(\delta),
\end{equation}
where
$\bm{\mathsf{R}}(\delta)=o(\delta^2)$; see Theorem \ref{SP_ZeroAvgForce_TransOnly}. Here, $\mathbb A$ is a symmetric positive definite matrix that depends only on the shape of $\mathscr B$, and is a measure of the drag exerted by $\mathscr L$ on $\mathscr B$ when the viscosity coefficient $\nu$ of $\mathscr L$ is very large. Furthermore, $\bm{G}$ is a spatial integral, involving the time average of the nonlinear terms evaluated at order $\delta$ and depending only on $\bm{F}$ and the physical parameters characterizing $\mathscr B$ and $\mathscr L$; see \eqref{SP_ZeroAvgForce_TransOnly} and Remark \ref{F0}.  Thus, at the order of $\delta^2$, $\bm{G}$ is the thrust responsible for the motion of $\mathscr B$. The evaluation of $\bm{G}$ is performed numerically in the particularly significant case where $\bm{F}$ results from the oscillation of the internal mass. To this end, with $a:={\rm diam}\,\mathscr B$, we define the quantity $h=a\sqrt{\pi/2\nu T}$ (Stokes number) and show that, after a suitable non-dimensionalization, $\bm{G}$ becomes a function of $h$ only. It turns out that, if $\mathscr B$ does not have fore-and-aft symmetry, like in the case of a drop-shaped body, then $\bm{G}\neq{\bf 0}$ for all allowed values of $h$, provided the oscillation of the mass is not too ``symmetric" in time; see Table \ref{tab1}. In such a case, $\bm{G}$ is an increasing, quadratic function of $h$ for large $h$, which means, in particular, that the speed eventually becomes a linear function of the frequency; see Table \ref{tab:frequency}. On the other hand, if $\mathscr B$ possesses fore-and-aft symmetry, $\bm{G}$ is calculated to be identically zero, for multiple choices of mass oscillation; see Table \ref{tab1}. This suggests that, for such bodies, propulsion may only take place at orders higher than $\delta^2$. We then analyze this question numerically, by a direct integration of the original equations. Our tests confirm the view above, in that they show that the same mass oscillations that produce $\bm{G}={\bf 0}$ now give a non-zero value for $\overline{\bm{\gamma}}$; see Tables \ref{tab1-nonlin} and \ref{tab45}, and Figure \ref{fig:ns-ellipse-23}. The question of whether this finding can also be proved analytically remains open and will be the subject of future work.

The plan of the paper is as follows. In Section \ref{sec:mathForm}, we give the mathematical formulation of the problem and show that propulsion is a nonlinear phenomenon. In Section \ref{ssec:fnSpcesand prelim}, we introduce the relevant function spaces and recall their main properties. In Section \ref{WkSolNonlnrProb}, we give the definition of weak solutions to our problem and show their existence in Theorem \ref{thm:SolLnrProb_TransOnly} for $\delta$ of restricted size. The following Section \ref{SCFP} is  devoted to find necessary and sufficient conditions for propulsion at the order $\delta^2$. Precisely, we show in Theorem  \ref{SP_ZeroAvgForce_TransOnly} that in the class of weak solutions, propulsion occurs at $O(\d^2)$ if and only if $\bm{G}\neq{\bf 0}$ and \eqref{th} holds. The final Section \ref{sec:NumericalResults} is entirely dedicated to the numerical analysis of the problem, where we find the results described above. 

\section{Formulation of the Problem}\label{sec:mathForm}
\begin{wrapfigure}{r}{0.37\linewidth}       
    \vspace{-12pt}
    \includegraphics[width=1\linewidth]{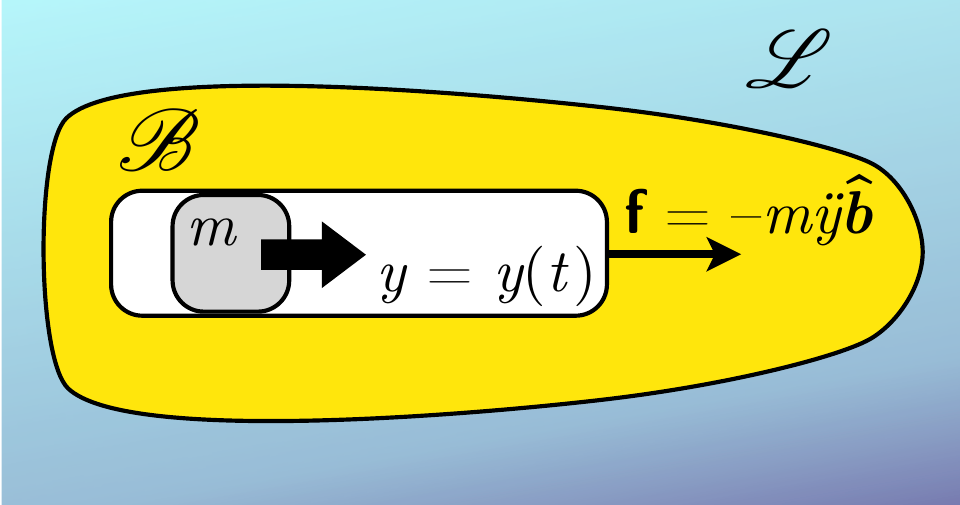}
    \caption{Schematic of the Driving Mechanism of $\B$.}
    \label{Schematic}
\end{wrapfigure}
    Consider a physical system, comprised of a rigid body $\B$ (closure of a bounded domain of $\mathbb R^3$ of class $C^2$) moving in a Navier-Stokes liquid $\L$ that fills the whole space outside $\B$, under the action of a force $\bm{\mathsf{f}}$. We assume that the motion of $\mathscr B$ is translatory, which can be achieved by applying on $\mathscr B$ a suitable torque preventing rotation. Concerning $\bm{\mathsf{f}}$, we suppose that it is time-periodic of period $T$ ($T$-periodic) and that its average, $\overline{\bm{\mathsf{f}}}$, over a period is zero:
\begin{equation}
\overline{\bm{\mathsf{f}}}:=\frac1T\int_0^T\bm{\mathsf{f}}(t)\,{\rm d}t={\bf 0}.
\label{zero}
\end{equation}
A remarkable special case that motivated our research is when $\bm{\mathsf{f}}$ is the result of the periodic displacement of a  mass $m$  in the interior of $\B$. More specifically, $\B$ houses a mechanism consisting of a cavity in which $m$ slides along a prescribed \textit{constant} (relative to $\mathscr B$) direction $\widehat{\bm{b}}$ in a time-periodic manner (see Figure \ref{Schematic}). As such, if $y=y(t)$ is the displacement of $m$ from some reference point, we know by Newton's law of motion, that this mechanism induces on $\mathscr B$ the force
\begin{equation}\label{eq:firsAppofForce}
    \bm{\mathsf{f}}(t)=-m\ddot{y}(t)\widehat{\bm{b}}\,.
\end{equation}
Thus, if $y(t)$ is $T$-periodic and sufficiently smooth, we deduce $\overline{\bm{\mathsf{f}}}={\bf 0}$.

Our analysis will be devoted to the general case, where $\bm{\mathsf{f}}$ is {\em any} (sufficiently smooth) $T$-periodic force satisfying \eqref{zero}. The equations of motion for the coupled body-liquid system $(\B,\L)$, referred to a frame attached to $\mathscr B$, are then given by \cite[Section 1]{galdi2002motion}:

\begin{equation}\label{eq:eom_coupled_system-TransOnly}
\begin{aligned}
    \left.\begin{array}{c}
        \displaystyle
        \frac{\partial\bm{v}}{\partial t}+(\bm{v}-\bm{\g})\cdot\nabla\bm{v}=\div\,\textbf{T}(\bm{v},p) \\
        \;\;\;\;\;\;\;\;\;\;\;\;\;\;\;\;\;\;\;\;\;\div\,\bm{v}=0 \\
    \end{array}\right\}&\;\;\;\;\;\;\;{\text{in}\;\Omega\times\mathbb{R}} \\\bm{v}=\bm{\g}&\;\;\;\;\;\;\;\text{on}\;\partial\Omega\times\mathbb{R} \\
    \lim_{|\bm{x}|\rightarrow\infty}\bm{v}(\bm{x},\cdot)=\textbf{0}\;&\;\;\;\;\;\;\;\text{in}\;\mathbb{R}. \\
        \begin{array}{c}
        \displaystyle
        \;\;M\dot{\bm{\g}}=\bm{f}-\int_{\partial\Omega}\textbf{T}(\bm{v},p)\cdot\bm{n}\,\text{d}S \\
    \end{array}&\;\;\;\;\;\;\;{\text{in}\;\mathbb{R}},
\end{aligned}
\end{equation}
where, $\bm{v}$ and $\rho p$ are the velocity and pressure fields of $\mathscr L$, respectively, with $\rho$ its density, whereas $\textbf{T}(\bm{v},p):=-p\textbf{1}+2\nu\textbf{D}(\bm{v})$, with \textbf{1} identity tensor,  $\nu$ kinematic viscosity and $\textbf{D}(\bm{v}):=\frac{1}{2}\left(\nabla\bm{v}+\left(\nabla\bm{v}\right)^{\top}\right)$, is the {Cauchy stress tensor}. Moreover, $\rho \,M$ represents the mass and $\bm{\gamma}$ the (translational) velocity of $\B$, and $\bm{n}$ stands for the outer unit normal to $\partial\Omega$. We also set
\begin{equation}\label{eq:scal_of_f_with_rho}
    \bm{\mathsf{f}}=:\rho \bm{f}\,.
\end{equation}
Thus, in view of \eqref{zero}, we have $\overline{\bm{f}}=\textbf{0}$.

As mentioned earlier, our goal is to find out when and how  $\bm{f}$ induces a propulsion on $\B$.  Since $\bm{f}$ is $T$-periodic, it is expected ---and shown in Theorem \ref{thm:ExistceWkSol_NonlnrProb}--- that problem \eqref{eq:eom_coupled_system-TransOnly} admits $T$-periodic solutions.  For such solutions, denote by $\bm{\eta}=\bm{\eta}(t)$ the position vector of the center of mass $G$ of $\mathscr B$ counted, say, from time $t=0$. We then have 
$$
    \bm{\eta}(t)=\int_0^t\bm{\gamma}(s)\,{\rm d}s+\bm{\eta}(0)\,.
$$
Because $\bm{\gamma}$ is $T$-periodic, it follows that the distance $d_T$ covered by $G$ in any interval of length $T$ is given by
$$
    d_T=\left|\int_0^T\bm{\gamma}(s)\,{\rm d}s\right|
\,.
$$
As a result, ensuring propulsion occurs is  equivalent to showing \eqref{gamma}; namely, that the time-average of $\bm{\gamma}$ over a period is non-zero.

\begin{remark}\label{ave} Our propulsion problem is strictly nonlinear. In fact, suppose $(\bm{v},p,\bm{\gamma})$ is a $T$-periodic solution to \eqref{eq:eom_coupled_system-TransOnly}. Denote by $\delta>0$ the magnitude of $\bm{f}$, so that $\bm{f}=\delta\,\bm{F}$. If we write $\bm{v}=\delta\,\bm{v}_0$, $p=\delta\,p_0$, $\bm{\gamma}=\delta\,\bm{\gamma}_0$, and disregard terms in $\delta^2$ and higher, we find that $(\bm{v}_0,p_0,\bm{\gamma}_0)$ obeys the following problem 
$$
\begin{aligned}
    \left.\begin{array}{c}
        \displaystyle
        \frac{\partial\bm{v}_0}{\partial t}=\div\,\textbf{T}(\bm{v}_0,p_0) \\
        \;\;\;\;\;\;\;\;\;\;\;\;\;\;\;\;\;\;\;\;\;\div\,\bm{v}_0=0 \\
    \end{array}\right\}&\;\;\;\;\;\;\;{\text{in}\;\Omega\times\mathbb{R}} \\\bm{v}_0=\bm{\g}_0&\;\;\;\;\;\;\;\text{on}\;\partial\Omega\times\mathbb{R} \\
    \lim_{|\bm{x}|\rightarrow\infty}\bm{v}_0(\bm{x},\cdot)=\textbf{0}\;&\;\;\;\;\;\;\;\text{in}\;\mathbb{R}. \\
        \begin{array}{c}
        \displaystyle
        \;\;M\dot{\bm{\g}}_0=\bm{F}-\int_{\partial\Omega}\textbf{T}(\bm{v}_0,p_0)\cdot\bm{n}\,\text{d}S \\
    \end{array}&\;\;\;\;\;\;\;{\text{in}\;\mathbb{R}},
\end{aligned}
$$  
Therefore, the averaged fields $(\overline{\bm{v}}_0,\overline{p}_0,\overline{\bm{\gamma}}_0)$ and  $\overline{\bm{F}}$ obey the boundary-value problem:
$$
\begin{aligned}
    \left.\begin{array}{c}
        \displaystyle
        \div\,\textbf{T}(\overline{\bm{v}}_0,\overline{p}_0)=\textbf{0} \\
        \div\,\overline{\bm{v}}_0=0 \\
    \end{array}\right\}&\;\;\;\;\;\;\;{\text{in}\;\Omega} \\
    \overline{\bm{v}}_0=\overline{\bm{\g}}_0&\;\;\;\;\;\;\;\text{on}\;\partial\Omega \\
    \lim_{|\bm{x}|\rightarrow\infty}\overline{\bm{v}}_0(\bm{x})=\textbf{0}\;& \\
        \begin{array}{c}
        \displaystyle
        \int_{\partial\Omega}\textbf{T}(\overline{\bm{v}}_0,\overline{p}_0)\cdot\bm{n}\,\text{d}S=\overline{\bm{F}} \\
    \end{array}&
\end{aligned}
$$
From classical results on the Stokes problem \cite[Section V.7]{galdi2011introduction}, we infer that, in a very large class of solutions, we may have $(\overline{\bm{v}_0},\nabla\overline{p_0},\overline{\bm{\gamma}_0})\neq({\bf 0},{\bf 0},{\bf 0})$ {\em if and only if} $\overline{\bm{F}}\neq{\bf 0}$; namely, $\overline{\bm{f}}\neq{\bf 0}$. Since, in our case $\overline{\bm{f}}={\bf 0}$, we thus conclude that propulsion can only occur at the order of $\delta^2$ or higher; that is, only when nonlinear effects are taken into account.
\end{remark}

\section{Function Spaces and Related Properties}\label{ssec:fnSpcesand prelim}

We begin to recall some basic notation. By $\Omega\subset\mathbb{R}^3$ we denote the complement in $\mathbb R^3$ of a compact domain $\mathscr B$ of class $C^2$ (the ``body"). We also set $R_*:=\text{diam}\,\mathscr{B}$ and take the origin of coordinates in the interior of $\mathscr B$. By $B_r$ we indicate the ball of radius $r>0$ in $\mathbb{R}^3$, and set
$$
    \Omega_R:=\Omega\cap B_R\,,\ \ R>R_*\,.
$$
For a domain $A\subseteq\mathbb{R}^3$,   $L^q(A)$  denotes the usual Lebesgue space endowed with the norm $\|\cdot\|_{L^q(A)}$ and, for $m\in\mathbb{N}$ and $q\in[1,\infty]$, $W^{m,q}(A)$ stands for the Sobolev space with norm $\|\cdot\|_{W^{m,q}(A)}$. We use $(\,\cdot\,,\,\cdot\,)_{L^2(A)}$ to denote the standard inner product in $L^2(A)$. Moreover, $D^{m,q}(A)$ will denote the homogeneous Sobolev space with semi-norm $|u|_{D^{m,q}(A)}:=\sum_{m=|\a|}\|D^{\a}u\|_{L^q(A)}$. Where the context is clear, we shall simply write $\|\,\cdot\,\|_q\equiv \|\,\cdot\,\|_{L^q(A)}$, $\|\,\cdot\,\|_{m,q}\equiv \|\,\cdot\,\|_{W^{m,q}(A)}$, $|\,\cdot\,|_{m,q}\equiv |\,\cdot\,|_{D^{m,q}(A)}$, and $(\,\cdot\,,\,\cdot\,)_2\equiv (\,\cdot\,,\,\cdot\,)_{L^2(A)}$. When any of the above function spaces are used with the subscript ``per", we shall mean that a function $u$ of this space has the additional property of being $T$-periodic; namely, $u(t+T)=u(t)$, for almost all $t\in\mathbb{R}$. Finally, given a function $w=w(t)$ defined in the interval $(0,T)$, we define its average:
$$
    \overline{w}:=\frac1T\int_0^Tw(t)\,\textrm{d}t.
$$

For $A\in \{\Omega,\Omega_R\}$ we introduce the set
\begin{equation*}
    \begin{aligned}
    \C(A)&:= \left\{\bm{\f}\in C_0^{\infty}(\overline{A}):
    \begin{array}{l}
    \text{$\div\,\bm{\f}=0$ in $A$;} \\
    \text{$\bm{\f}=\bm{\g}_{\bm{\f}}$\  in a neighborhood of $\partial A$, for some $\bm{\g}_{\bm{\f}} \in \R^3$;} \\
    \text{$\bm{\f}=\bm{0}$ in a neighborhood of $\partial B_R$ if $A\equiv \Omega_R$}
    \end{array}
    \right\},
    \end{aligned}
\end{equation*}
so that, considering the inner products,
$$
    \left(\u,\w\right)_{\K(A)}:=\left( \u,\w\right)_{L^2(A)}+M\bm{\g}_{\u}\cdot\bm{\g}_{\w},\;\;\;\text{for all}\;\u,\w\in \C(A),
$$
and
$$
    \left(\u,\w\right)_{\H(A)}:=\int_{A}\textbf{D}(\u):\textbf{D}(\w),\;\;\;\text{for all}\;\u,\w\in \C(A),
$$
with associated norms,
$$
    \|\u\|_{\K(A)}:=\left(\|\u\|^2_{L^2(A)}+M|\bm{\g}_{\u}|^2\right)^{1/2},\;\;\;\text{for all}\;\u\in \C(A),
$$
and
$$
    \|\u\|_{\H(A)}:=\|\textbf{D}(\u)\|_{L^2(A)},\;\;\;\text{for all}\;\u\in \C(A),
$$
respectively, we can define the completions,
$$
    \K(\O_R):=\overline{\C(\O_R)}^{\|\cdot\|_{\K(\O_R)}}\;\;\;\;\;\text{and}\;\;\;\;\;\H(A):=\overline{\C(A)}^{\|\cdot\|_{\H(A)}}.
$$
We shall often simply use the subscripts $\K$ and $\H$ whenever $A\equiv\O$.

It can be shown (see   \cite[Lemma 4.11]{galdi2002motion}), that
\begin{equation*}
    \begin{aligned}
    \H(\O)&:= \left\{\v\in W^{1,2}_\text{loc}(\overline{\O})\cap L^6(\O):
    \begin{array}{l}
    \text{$\textbf{D}(\v)\in L^2(\O)$, $\div\,\v=0$, and} \\
    \text{$\v=\bm{\g}_{\v}$, for some $\bm{\g}_{\bm{\v}} \in \R^3$ on $\partial \O$}
    \end{array}
    \right\}.
    \end{aligned}
\end{equation*}
Likewise, for the ``local spaces", the following characterizations hold:
\begin{equation*}
\begin{aligned}
    \K(\O_R) &= \{\u\in L^2(\O_R): \text{$\div\,\u=0$ in $\O_R$; $\u\cdot\bm{n}=\bm{\g}_{\u}\cdot\bm{n}$ around $\partial\O$; $\u\cdot\bm{n}=0$ around $\partial B_R$}\}; \\
    \H(\Omega_R)  &= \{\u\in W^{1,2}(\O_R): \text{$\div\,\u=0$ in $\O_R$; $\u=\bm{\g}_{\u}$ around $\partial\O$; $\u=\bm{0}$ around $\partial B_R$}\}.
\end{aligned}
\end{equation*}
It is known (and easy to check) that $\K(\O_R)$ and $\H(A)$ are Hilbert spaces when endowed with scalar products  $(\cdot,\cdot)_{\K(\O_R)}$ and $(\cdot,\cdot)_{\H(A)}$, respectively. For $m\in\mathbb{N}\cup\{\infty\}$ and fixed $T>0$, we introduce the test function spaces
\begin{equation*}
    \begin{aligned}
    \C^m_\text{per}(A\times\R)&:= \left\{\bm{\f}\in C^m(A\times\R):
    \begin{array}{l}
    \text{$\div\,\bm{\f}=0$ in $A$; $\bm{\f}$ is $T$-periodic;} \\
    \text{$\bm{\f}(\x,\cdot)=\bm{\g}_{\bm{\f}}(\cdot)$, for some $\bm{\g}_{\bm{\f}} \in C_\text{per}^m(\R)$,} \\
    \text{for all $\x$ in a neighbourhood of $\partial \O$;} \\
    \text{there exists $r>\text{diam}\,\B$, such that $\bm{\f}(\x,t)=0$,} \\
    \text{for all $\x\in \overline{B}_{r}^c$ and all $t\in\R$, where $r<R$ if $A\equiv \O_R$}
    \end{array}
    \right\}
    \end{aligned}
\end{equation*}
where we use $\C^m_\text{per}(A\times[0,T])$ to denote the functions of $\C^m_\text{per}(A\times\mathbb{R})$ restricted to $[0,T]$. Similarly, we will use $C^m_\text{per}([0,T])$ to denote the functions of $C^m_\text{per}(\mathbb{R})$ restricted to $[0,T]$.

Given $q\in[1,\infty]$, we establish the conventions
$$
    \widehat{L}^q_\text{per}(\R):=\{\bm{\xi}\in L^q_\text{per}(\R):\overline{\bm{\xi}}=\textbf{0}\}\;\;\;\;\;\text{and}\;\;\;\;\;\widehat{W}_\text{per}^{1,q}(\R):=\{\bm{\xi}\in \widehat{L}^q_\text{per}(\R):\dot{\bm{\xi}}\in L^q_\text{per}(\R)\}
$$
and, using the shorthand $X(Y)\equiv X(\R;Y)$ for function spaces $X$ and $Y$ (e.g. $W^{m,q}_\text{per}(L^r(A))\equiv W^{m,q}_\text{per}(\R;L^r(A))$, dropping the domain $A$ where the context is clear, we define the space
$$
    \mathcal{W}:=\left[\widehat{W}^{1,2}_\text{per}(L^2)\cap \widehat{L}^2_\text{per}(W^{2,2}\cap\H)\right]
$$
with norm
$$
    \|\bm{u}\|_{\mathcal{W}}:=\|\bm{u}\|_{W^{1,2}(L^2)}+\|\bm{u}\|_{L^2(W^{2,2})}.
$$
With regards to the former space, we have the following embedding result \cite[Theorem 2.1]{solonnikov1973EstofSol}.

\begin{lemma}\label{thm:SolonnikovEmbPerZAFns}
    Let $q\in[1,\infty)$. Then,  the following embedding holds for all $r,s\in [q,\infty]$:
    $$
        \mathcal{W}\hookrightarrow L^r(L^s),\;\;\;\;\;\;\;\frac{3}{s}+\frac{2}{r}>\frac{1}{2}.
    $$
\end{lemma}

We conclude this section with the following lemma, containing a collection of important estimates pertaining to the spaces $\H(A)$ (see \cite[Section 4]{galdi2002motion}).

\begin{lemma}\label{thm:KornsIdentity_coupledSyst}
    For $R>R_*$, let $A\in\{\Omega,\Omega_R\}$ and $\bm{u}\in\H(A)$. Then
    \begin{equation}\label{eq:KornsIdentity_coupledSyst}
        \|\nabla \bm{u}\|_{L^2(A)}=\sqrt{2}\|\bm{u}\|_{\H(A)}
    \end{equation}
    and there exist $c_1,c_2>0$, independent of $A$, such that, for all $\bm{u}\in\H(A)$, the following inequalities hold:
    \begin{eqnarray}
        |\bm{\g}_{\bm{u}}|&\leq& c_1\|\bm{u}\|_{\H(A)}\label{eq:estOnGamma_u} \\
        \|\bm{u}\|_{L^6(A)} &\leq& c_2 \|\bm{u}\|_{\H(A)}.\label{eq:sobolevInHO_R}
    \end{eqnarray}
\end{lemma}

\section{\textbf{Weak Solutions to a Nonlinear Problem}}\label{WkSolNonlnrProb}

To address the question of propulsion, we will apply a slightly similar, though much simpler, procedure to the one used in \cite[Section 6]{galdi2020self}. For this purpose, set
$$
    \delta:= \|\bm{f}\|_{L^2(0,T)}\,,\ \ \bm{F}:=\bm{f}/\|\bm{f}\|_{L^2(0,T)}
$$
and consider the following linearized problem,
\begin{equation}\label{eq:eom_coupledSyst_fEQ0Case-1}
\begin{aligned}
    \left.\begin{array}{c}
        \displaystyle
        \;\;\;\;\;\;\;\;\;\;\;\;\;\;\;\;\;\;\;\;\,\frac{\partial\bm{V}_0}{\partial t}=\div\textbf{T}(\bm{V}_0,p_0) \\
        \div\,\bm{V}_0=0 \\
    \end{array}\right\}&\;\;\;\;\;\;\;{\text{in}\;\O\times\R} \\
    \bm{V}_0(\bm{x},t)=\bm{\xi}_0(t)&\;\;\;\;\;\;\;\text{on}\;\partial\O\times\R \\
    M\dot{\bm{\xi}}_0=\bm{F}-\int_{\partial\O} \textbf{T}(\bm{V}_0,p_0)\cdot\bm{n}\;\text{d}S&\;\;\;\;\;\;\;\text{in}\;\R,
\end{aligned}
\end{equation}
which is known to have a solution. In fact, we have the following \cite[Lemma 5.2]{galdi2020self}.

\begin{thm}\label{thm:SolLnrProb_TransOnly}
    For any $\bm{F}\in \widehat{L}_\text{per}^2(\R)$, there exists a unique solution $(\bm{V}_0,p_0,\bm{\xi}_0)\in \mathcal{W}\times \widehat{L}_\text{per}^2(D^{1,2})\times \widehat{W}_\text{per}^{1,2}(\R)$ to problem (\ref{eq:eom_coupledSyst_fEQ0Case-1}). In addition, this solution satisfies the estimate~\footnote{Observe that $\|\bm{F}\|_{L^2(0,T)}=1$.}
    $$
        \|\bm{V}_0\|_{\mathcal{W}}+\|p_0\|_{L^2(D^{1,2})}+\|\bm{\xi}_0\|_{W^{1,2}(0,T)}\leq C_0\,,
    $$
with $C_0>0$ depending only on $\Omega$ and $T$.
\end{thm}
We shall then look for a solution to $(\ref{eq:eom_coupled_system-TransOnly})$ ``around" the solution $(\bm{V}_0,p_0,\bm{\xi}_0)$. Precisely, for each $\d>0$, we introduce the  formal decompositions of the fields $(\v,p,\bm{\g})$ in (\ref{eq:eom_coupled_system-TransOnly}):
\begin{equation}\label{eq:decomp_1stappearance}
    \v:=\u+\d\bm{V}_0,\;\;\;\;\;\;\;\bm{\g}:=\bm{\xi}+\d\bm{\xi}_0,\;\;\;\;\;\;\;p:=\pi+\d p_0.
\end{equation}

Thus, substituting (\ref{eq:decomp_1stappearance}) into (\ref{eq:eom_coupled_system-TransOnly}), we obtain the following nonlinear problem in the unknowns $(\u,\bm{\xi},\pi)$:
\begin{equation}\label{eq:eom_coupledSyst_fEQ0Case-2}
\begin{aligned}
    \left.\begin{array}{c}\medskip
        \displaystyle
        \frac{\partial\u}{\partial t}+(\u-\bm{\xi})\cdot\nabla\u+\d\left[ (\u-\bm{\xi})\cdot\nabla\bm{V}_0+ (\bm{V}_0-\bm{\xi}_0)\cdot\nabla\u\right]\\ \medskip
        \hspace*{8cm}=\div\textbf{T}(\u,\pi)+\d^2\bm{F}_0 \\
        \hspace*{4.2cm}\div\,\u=0 \\
    \end{array}\right\}&\;\;\;\;\;\;\;{\text{in}\;\O\times\R} \\
    \u(\bm{x},t)=\bm{\xi}(t)&\;\;\;\;\;\;\;\text{on}\;\partial\O\times\R \\
    M\dot{\bm{\xi}}=-\int_{\partial\O} \textbf{T}(\u,\pi)\cdot\bm{n}\;\text{d}S&\;\;\;\;\;\;\;\text{in}\;\R,
\end{aligned}
\end{equation}
where 
$$
\bm{F}_0:=-(\bm{V}_0-\bm{\xi}_0)\cdot\nabla\bm{V}_0\,.
$$
\begin{remark}\label{F0} We would like to emphasize that, in problem \eqref{eq:eom_coupledSyst_fEQ0Case-2}, $\bm{F}_0$ is the only driving data, in the sense that if $\bm{F}_0\equiv{\bf 0}$, then $\bm{u}\equiv\nabla\pi\equiv\bm{\xi}\equiv{\bf 0}$ is a solution. Moreover, $\bm{F}_0$ is a genuine nonlinear contribution, since it is just the nonlinearity in \eqref{eq:eom_coupled_system-TransOnly}$_1$  evaluated along the solution $(\bm{V}_0,\bm{\xi}_0)$ to \eqref{eq:eom_coupledSyst_fEQ0Case-1}.
\end{remark}
The aim of this section will be to prove existence of solutions to a suitable weak formulation of  problem  \eqref{eq:eom_coupledSyst_fEQ0Case-2}. To this end, we first dot-multiply (\ref{eq:eom_coupledSyst_fEQ0Case-2}$)_1$ by arbitrary $\bm{\f}\in \C_\text{per}^1(\Omega\times\mathbb{R})$ and integrate by parts using (\ref{eq:eom_coupledSyst_fEQ0Case-2}$)_\text{2,3,4}$ and also periodicity to get
\begin{equation}\label{eq:wkForm_coupledSyst_fEQ0Case-2}
    \begin{aligned}
        &\int_0^T {\Bigg[}\left(\u,\frac{\partial\bm{\f}}{\partial t}\right)_{\K} -\d\left[ \left((\u-\bm{\xi})\cdot\nabla\bm{V}_0,\bm{\f}\right)_{2}+ \left((\bm{V}_0-\bm{\xi}_0)\cdot\nabla\u,\bm{\f}\right)_{2}\right] \\
        &\;\;\;\;\;\;\;\;\;\;\;\;\;\;\;\;\;\;\;\;\;\;\;\;\;\;\;\;\;\;\;\;\;\;\;\;\;\;\;\;\;\;\;\;\;\;\;\;\;- \left((\u-\bm{\xi})\cdot\nabla\u,\bm{\f}\right)_{2}-2\nu\left( \u,\bm{\f}\right)_{\H}+ \d^2(\bm{F}_0,\bm{\f})_{2}{\Bigg]}\text{d}t=0,
    \end{aligned}
    \end{equation}

With this in hand, we can now introduce the following definition.

\begin{definition}\label{def:wkSol_NonLrnProbinDelta}
    Let $\d>0$ and $(\bm{V}_0,\bm{\xi}_0)\in \mathcal{W}\times \widehat{W}_\text{per}(\R)$. Then, the pair $(\u, \bm{\xi})$ is said to be a {\em $T$-periodic weak solution} to problem (\ref{eq:eom_coupledSyst_fEQ0Case-2}) if
\begin{enumerate}
    \item $\u\in L_\text{per}^2(\mathbb{R};\H(\Omega))$ and $\bm{\xi}\in L^2_\text{per}(\mathbb{R})$ with $\u=\bm{\xi}$ on $\partial\O$ (in the trace sense)\footnote{{In fact, due to (\ref{eq:estOnGamma_u}), the condition $\bm{\xi}\in L^2_\text{per}(\mathbb{R})$ is automatic if we take $\bm{\xi}\equiv\bm{\xi}_{\bm{u}}$.}};
    
    \item $(\u, \bm{\xi})$ verifies (\ref{eq:wkForm_coupledSyst_fEQ0Case-2}) for all $\bm{\f}\in\C^1_\text{per}(\Omega\times \mathbb{R})$.
\end{enumerate}
\end{definition}

To prove the existence of $T$-periodic weak solutions to (\ref{eq:eom_coupledSyst_fEQ0Case-2}), we first prove existence of a ``localized solution" on a bounded domain $\O_R$, for any $R>R_*$ using the Galerkin method and then use the so-called ``invading domains" technique to pass to the limit in $R$ and obtain a solution on the whole domain $\O$. For the former step, we shall require the following special basis \cite[Lemma 3.1]{galdi2009motion}.\looseness=-1

\begin{lemma}\label{thm:onBasis_coupledSyst}
    Let $R>R_*$. Then, there is an orthonormal basis $\left\{\bm{\Psi}_k\right\}_{k\in\N}\subset \C(\O_R)$ of $\K(\O_R)$ such that, if
    $$
        \mathcal{X}_N(\O_R\times\R):=\left\{\sum_{k=1}^N\chi_k(t)\bm{\Psi}_k(\x) : \chi_k\in C_{0,\text{per}}^1(\R) \right\},
    $$
    then given any $\bm{\f}\in\C_{\text{per}}^1(\O_R\times\R)$, for all $\e>0$, there exists $N\in\N$ and a function $\bm{\f}_N\in\mathcal{X}_N(\O_R\times\R)$ satisfying the following:
    \begin{multicols}{2}
    \begin{enumerate}
        \item $\displaystyle \max_{[0,T]}\left\|\bm{\f}_N(t)-\bm{\f}(t)\right\|_{C^2(\O_R)}<\e$
        
        \item $\displaystyle \max_{[0,T]}\left\|\frac{\partial\bm{\f}_N}{\partial t}(t)-\frac{\partial\bm{\f}}{\partial t}(t)\right\|_{C^0(\O_R)}<\e.$
        
        \item $\displaystyle\max_{[0,T]}\left|\bm{\g}_{\bm{\f}_N}-\bm{\g}_{\bm{\f}}\right|<\e$
        
        \item $\displaystyle\max_{[0,T]}\left|\dot{\bm{\g}}_{\bm{\f}_N}-\dot{\bm{\g}}_{\bm{\f}}\right|<\e$
    \end{enumerate}
    \end{multicols}
\end{lemma}

\begin{thm}\label{thm:ExistceWkSol_NonlnrProb}
    Let $(\bm{V}_0,\bm{\xi}_0)\in \mathcal{W}\times \widehat{W}_\text{per}^{1,2}(\R)$. Then, there exists $\d_0>0$ such that, for all $\d\in (0,\d_0)$, there is at least one $T$-periodic corresponding weak solution $(\u,\bm{\xi})$ to problem (\ref{eq:eom_coupledSyst_fEQ0Case-2}) satisfying the estimate
    \begin{equation}\label{eq:UnifEst_WkSol_NonlnrProb_Omega}
        \|\u\|_{L^2(\H)}+\|\bm{\xi}\|_{L^2(0,T)}\leq c\,\d^2.
    \end{equation}
\end{thm}
\begin{proof}
    As mentioned earlier, the proof shall be accomplished in two principal steps, using a standard approach (for details, we refer the reader to \cite[Section 3]{galdi2009motion}, \cite[Theorem 1]{galdi2018time}). First, given any $R>R_*$, using the Galerkin method, we prove, for $\d$ sufficiently small, existence of a pair $(\u_R,\bm{\xi}_R)\in L_\text{per}^2(\mathbb{R};\H(\O_R))\times L^2_\text{per}(\mathbb{R})$ with $\u_R=\bm{\xi}_R$, satisfying
    \begin{equation}\label{eq:LOCALwkForm_coupledSyst_fEQ0Case-2}
    \begin{aligned}
        &\int_0^T {\Bigg[}\left(\u_R,\frac{\partial\bm{\f}}{\partial t}\right)_{\K(\O_R)} -\d\left[ \left((\u_R-\bm{\xi}_R)\cdot\nabla\bm{V}_0,\bm{\f}\right)_{L^2(\O_R)}+ \left((\bm{V}_0-\bm{\xi}_0)\cdot\nabla\u_R,\bm{\f}\right)_{L^2(\O_R)}\right] \\
        &\;\;\;\;\;\;\;\;\;\;\;\;\;\;\;\;\;\;\;\;\;\;\;\;\;- \left((\u_R-\bm{\xi}_R)\cdot\nabla\u_R,\bm{\f}\right)_{L^2(\O_R)}-2\nu\left( \u_R,\bm{\f}\right)_{\H(\O_R)}+ \d^2(\bm{F}_0,\bm{\f})_{L^2(\O_R)}{\Bigg]}\text{d}t=0,
    \end{aligned}
    \end{equation}
    for all $\bm{\f}\in\C^1_\text{per}(\Omega_R\times \mathbb{R})$ as well as the estimate
    \begin{equation}\label{eq:NonlnrProb_UnifEst_inOmegaR}
        \|\u_R\|_{L^2(0,T;\H(\O_R))}+\|\bm{\xi}_R\|_{L^2(0,T)}\leq c\d^2,
    \end{equation}
    for some constant $c>0$ independent of $R$. Then, in the second step, we show that these solutions converge, in a suitable sense, to a solution $(\u,\bm{\xi})$ on the whole domain $\O$.

    \hfill

    \textit{Step 1.} With $\{\bm{\Psi}_k\}_{k\in\N}\subset \C(\O_R)$, the basis for $\K(\O_R)$ given by Lemma \ref{thm:onBasis_coupledSyst}, for $\d$ sufficiently small, we begin by looking for ``approximating solutions" of the form
    \begin{equation}\label{eq:approxWkSol_NonlnrProb}
        \u_n=\sum_{k=1}^n\a_{nk}(t)\bm{\Psi}_k(\x),\;\;\;\;\;\;\;\bm{\xi}_n=\sum_{k=1}^n\a_{nk}(t)\bm{\xi}_{\bm{\Psi}_k},
    \end{equation}
    satisfying
    \begin{equation}\label{eq:approxWkForm_NonlnrProb}
    \begin{aligned}
        &\left(\frac{\partial\u_n}{\partial t},\bm{\Psi}_k\right)_{\K(\O_R)} +\d\left[ \left((\u_n-\bm{\xi}_n)\cdot\nabla\bm{V}_0,\bm{\Psi}_k\right)_{L^2(\O_R)}+ \left((\bm{V}_0-\bm{\xi}_0)\cdot\nabla\u_n,\bm{\Psi}_k\right)_{L^2(\O_R)}\right] \\
        &\hspace{2cm}+ \left((\u_n-\bm{\xi}_n)\cdot\nabla\u_n,\bm{\Psi}_k\right)_{L^2(\O_R)}+2\nu\left( \u_n,\bm{\Psi}_k\right)_{\H(\O_R)}- \d^2(\bm{F}_0,\bm{\Psi}_k)_{L^2(\O_R)}=0,
    \end{aligned}
    \end{equation}
and the estimates
\begin{equation}\label{eq:PfNonLnrProbunif_bd_coupledSyst}
    \|\u_n\|_{L^2(0,T;\H(\O_R))}\leq C_1\d^2\;\;\;\;\;\text{and}\;\;\;\;\;\|\u_n\|_{L^{\infty}(0,T;\K(\O_R))}\leq C_2,
\end{equation}
where the constant $C_1>0$ is independent of both $R$ and $n$ and $C_2>0$ depends only on $R$. To this end, substituting (\ref{eq:approxWkSol_NonlnrProb}) into (\ref{eq:approxWkForm_NonlnrProb}), using the orthonormality of the basis vectors, yields a system of ordinary differential equations, which, by the standard theory of ordinary differential equations, admits a unique solution $\bm{\a}_n(t):=(\a_{n1}(t),...,\a_{nn}(t))\in W^{1,2}(0,T_n)$ for some $T_n>0$, where we can take $T_n=T$, provided $\bm{\a}_n$ is uniformly bounded. We now show that, for sufficiently small $\d$, this is indeed the case.

Multiplying (\ref{eq:approxWkForm_NonlnrProb}) by $\a_{nj}$ and summing over $j$, we obtain
\begin{equation}\label{eq:PfExistNonlnrProb_aPrioriEst-UnInt}
    \frac{d}{dt}\|\u_n\|_{\K(\O_R)}^2+4\nu\|\u_n\|_{\H(\O_R)}^2+2\d\left((\u_n-\bm{\xi}_n)\cdot\nabla\bm{V}_0,\u_n\right)_{L^2(\O_R)}=2\d^2(\bm{F}_0,\u_n)_{L^2(\O_R)}.
\end{equation}

Now, applying the H\"older and Young Inequalities along with Lemma \ref{thm:KornsIdentity_coupledSyst}, we show the following two estimates:
\begin{equation}\label{eq:PfExistNonlnrProb_SolonnEmbApplied2NonlnrTerms}
\begin{aligned}
    \left|\left((\u_n-\bm{\xi}_n)\cdot\nabla\bm{V}_0,\u_n\right)_{L^2(\O_R)}\right| &\leq c_1\left(\|\bm{V}_0\|_{L^2(\O_R)}+\|\bm{V}_0\|_{L^3(\O_R)}\right)\|\u_n\|_{\H(\O_R)}^2 \\
    \left|(\bm{F}_0,\u_n)_{L^2(\O_R)}\right| &\leq c_2\left(\|\bm{V}_0\|_{L^4(\O_R)}+|\bm{\xi}_0|\,\|\bm{V}_0\|_{L^2(\O_R)}\right)\|\u_n\|_{\H(\O_R)},
\end{aligned}
\end{equation}
for constants $c_1,c_2>0$ independent of $R$. We next observe that, by Lemma \ref{thm:SolonnikovEmbPerZAFns}, and the embedding $W^{1,2}(0,T)\subset L^\infty(0,T)$, we deduce
\begin{equation}\label{GP}\begin{array}{ll}\medskip
\|\bm{V}_0\|_{L^\infty(L^q)}\le c_3\,\|\bm{V}_0\|_{\mathcal W}\le c_3\,C_0\,,\ \ \mbox{for all $q\in [2,6)$}\,,
\\
\|\bm{\xi}_0\|_{L^\infty(0,T)}\le c_4\|\bm{\xi}_0\|_{L^\infty(0,T)}\le c_4\,C_0
\end{array}
\end{equation}
where $C_0$ is the constant entering the estimate in Theorem \ref{thm:SolLnrProb_TransOnly}. Therefore, combining \eqref{eq:PfExistNonlnrProb_aPrioriEst-UnInt} and \eqref{eq:PfExistNonlnrProb_SolonnEmbApplied2NonlnrTerms} with \eqref{GP} and taking $\delta\in (0,\delta_0)$, for a suitable $\delta_0>0$, we get
\begin{equation}\label{delta}
 \frac{d}{dt}\|\u_n\|_{\K(\O_R)}^2+c_5\,\|\u_n\|_{\H(\O_R)}^2\le c_6\,\delta^2\,,
\end{equation}
where $c_5$ and $c_6$ are independent of $R$. Integrating from 0 to $t$, this differential inequality, we deduce, in particular,
\begin{equation*}
    |\bm{\a}_n(t)|=\|\u_n(t)\|_{\K(\O_R)}\leq c_6\,T\,\d^2+|\bm{a}|^2,\;\;\;\;\;\;\;\text{for all $\d\in(0,\d_0)$ and all $t\in(0,T]$}\,,
\end{equation*}
where $\bm{a}=\bm{\a}_n(0)$ are given initial data. Hence, we can take $T_n=T$.

We now proceed to show that the coefficients $\a_{n1},...,\a_{nn}$ can be taken to be $T$-periodic. Employing Poincar\'e inequality combined with (\ref{eq:estOnGamma_u}) on the left-hand side of \eqref{delta}, we infer with $\kappa=\kappa(R)>0$
\begin{equation}\label{eq:PfNonlnrProb_DiffIne}
     \frac{d}{dt}\|\u_n\|_{\K(\O_R)}^2+\kappa\|\u_n\|_{\K(\O_R)}^2 \leq c_7,\;\;\;\;\;\;\;c_7=c_7(\d_0)>0,
\end{equation}
which, by the classical Gr\"onwall inequality, implies
\begin{equation}\label{eq:PfNonlnrProb_DiffIne_Intd}
    |\bm{\a}_n(T)|^2=\|\u_n(T)\|_{\K(\O_R)}^2\leq e^{-\kappa T}|\bm{a}|^2+\frac{c_7}{\kappa}.
\end{equation}
Now, by uniqueness of the solution $\bm{\a}_n$, the flow map $F:\R^n\rw\R^n$, given by $F(\bm{a}):=\bm{\a}_n(T)$, is well-defined and continuous. Moreover, taking $r^2 \geq c_7/\kappa(1-e^{-\kappa T})$, one easily finds, using (\ref{eq:PfNonlnrProb_DiffIne_Intd}), that $F$ maps the closed ball $\overline{B_r}$ to itself. As such, by Brouwer's fixed-point theorem, $F$ has a fixed point, say $\widetilde{\bm{a}}\in B_r$, so that $\bm{\a}_n(T)=F(\widetilde{\bm{a}})=\widetilde{\bm{a}}=\bm{\a}_n(0)$, making $\bm{\a}_n$ $T$-periodic for this choice of initial value. Therefore, we conclude that we can take $\u_n$ and $\bm{\xi}_n$ to be $T$-periodic.

Let us now establish the estimates in (\ref{eq:PfNonLnrProbunif_bd_coupledSyst}). The first one follows at once from \eqref{delta}, after integrating both sides from 0 to $T$, and using the periodicity of $\u_n$. We therefore turn our attention to (\ref{eq:PfNonLnrProbunif_bd_coupledSyst}$)_2$. First, by the Mean Value Theorem, we find a number $\widetilde{T}\in(0,T)$, such that
$$
    \frac{1}{T}\int_0^T\|\u_n(t)\|^2_{\H(\O_R)}\,\text{d}t=\|\u_n(\widetilde{T})\|^2_{\H(\O_R)}.
$$
Using this result combined with the Poincar\'e inequality, (\ref{eq:estOnGamma_u}), and (\ref{eq:PfNonLnrProbunif_bd_coupledSyst}$)_1$, we obtain
\begin{equation}\label{eq:prelimEst-on-u_n_coupledSyst}
    \|\u_n(\widetilde{T})\|_{\K(\O_R)}\leq c_8,\;\;\;\;\;\;\;c_8=c_8(R,\d_0)>0.
\end{equation}
Then (\ref{eq:PfNonLnrProbunif_bd_coupledSyst}$)_2$, follows by integrating \eqref{delta} from $\widetilde{T}$ to $t>\widetilde{T}$, and applying (\ref{eq:prelimEst-on-u_n_coupledSyst}).

Now, consider the sequences $\{\u_n\}_{n\in\N}$ and $\{\bm{\xi}_n\}_{n\in\N}$. Using (\ref{eq:PfNonLnrProbunif_bd_coupledSyst}) combined with standard procedures (see, for instance, \cite[Section 3]{galdi2000introduction}), one proves the existence of a pair $(\u_R,\bm{\xi}_R)$ of vector fields, such that $\bm{\u}_R=\bm{\xi}_R$ on $\partial\O$ and satisfying (up to subsequence) the convergence properties
\begin{enumerate}
    \item $\u_n \xrightharpoonup{\;\;\;\;\;} \u_R$ in $L^2_\text{per}(\R;\H(\O_R))$,
    \item $\u_n\longrightarrow \u_R$ in $L_\text{per}^2(\R;\K(\Omega_R))$,  and
    \item $\bm{\xi}_n\longrightarrow \bm{\xi}_R$ in $L_\text{per}^2(\R)$,
\end{enumerate}
as well as the estimate (\ref{eq:NonlnrProb_UnifEst_inOmegaR}). 

Finally, choosing any $N\in\N$, after multiplying (\ref{eq:approxWkForm_NonlnrProb}) by arbitrary $\chi_{j}\in C_{0,\text{per}}^1(\R)$, summing over $j=1,...,N$, and then integrating from $0$ to $T$, we immediately conclude that $(\u_n,\bm{\xi}_n)$ satisfies the following ``approximating weak form"
\begin{equation}\label{eq:LOCALwkForm_coupledSyst_fEQ0Case-2MultThroughbyphiN}
    \begin{aligned}
        &\int_0^T {\Bigg[}\left(\u_n,\frac{\partial\bm{\f}_N}{\partial t}\right)_{\K(\O_R)} -\d\left[ \left((\u_n-\bm{\xi}_n)\cdot\nabla\bm{V}_0,\bm{\f}_N\right)_{L^2(\O_R)}+ \left((\bm{V}_0-\bm{\xi}_0)\cdot\nabla\u_n,\bm{\f}_N\right)_{L^2(\O_R)}\right] \\
        &\;\;\;\;\;\;\;\;\;\;\;\;\;\;\;\;\;- \left((\u_n-\bm{\xi}_n)\cdot\nabla\u_n,\bm{\f}_N\right)_{L^2(\O_R)}-2\nu\left( \u_n,\bm{\f}_N\right)_{\H(\O_R)}+ \d^2(\bm{F}_0,\bm{\f}_N)_{L^2(\O_R)}{\Bigg]}\text{d}t=0,
    \end{aligned}
    \end{equation}
for all $\bm{\f}_N\in \mathcal{X}_N(\O_R\times\R)$. Then, employing similar procedures as those used in \cite[Section 3]{galdi2009motion}, \cite[Theorem 1]{galdi2018time}, we use the convergence properties (i)-(iii) as well as those given in Lemma \ref{thm:onBasis_coupledSyst} to pass to the limit in (\ref{eq:LOCALwkForm_coupledSyst_fEQ0Case-2MultThroughbyphiN}), first as $n\rw\infty$ and then as $N\rw\infty$.

\hfill

\textit{Step 2.} Now, take a sequence of increasing radii $\{R_k\}_{k\in\N}$, $R_1>R_*$, $\lim_{k\rw\infty}R_k=\infty$, and define the shorthand
$$
    \u_k:=\u_{R_k},\;\;\;\bm{\xi}_k:=\bm{\xi}_{R_k},\;\;\;\O_k:=\O_{R_k}.
$$
In this notation, by Step 1, for each $k\in\N$, we have a solution $(\u_k,\bm{\xi}_k)$ on the domain $\O_k$. So now, fix $\bm{\f}\in \C_\text{per}^1(\O\times\R)$ and take $k_0\in\N$ sufficiently large, so that $\text{supp}\,\bm{\f}\subset\O_{k_0}=:\O_0$. Then, thanks to Step 1, after extending $\u_k$ to be zero outside $\overline{B}_{R_k}$, we have
\begin{equation}\label{eq:PfNonlnrProb_UnifBd4InvDom_inK}
    \|\u_k\|_{L^2(0,T;\H(\O_R))}+\|\bm{\xi}_k\|_{L^2(0,T)}\leq c\,\d^2,
\end{equation}
and also that $(\u_k,\bm{\xi}_k)$ satisfies
\begin{equation}\label{eq:LOCALwkForm_coupledSyst_fEQ0Case-finalApprox}
    \begin{aligned}
        &\int_0^T {\Bigg[}\left(\u_k,\frac{\partial\bm{\f}}{\partial t}\right)_{\K(\O)} -\d\left[ \left((\u_k-\bm{\xi}_k)\cdot\nabla\bm{V}_0,\bm{\f}\right)_{L^2(\O)}+ \left((\bm{V}_0-\bm{\xi}_0)\cdot\nabla\u_k,\bm{\f}\right)_{L^2(\O)}\right] \\
        &\;\;\;\;\;\;\;\;\;\;\;\;\;\;\;\;\;- \left((\u_k-\bm{\xi}_k)\cdot\nabla\u_k,\bm{\f}\right)_{L^2(\O)}-2\nu\left( \u_k,\bm{\f}\right)_{\H(\O)}+ \d^2(\bm{F}_0,\bm{\f})_{L^2(\O)}{\Bigg]}\text{d}t=0,
    \end{aligned}
    \end{equation}
for all $k \geq k_0$.

Now, by the weak compactness property of the space $L^2_\text{per}(\R;\H(\O))$, from (\ref{eq:PfNonlnrProb_UnifBd4InvDom_inK}), we extract $(\u,\bm{\xi})\in L^2_\text{per}(\R;\H(\O))\times L^2_\text{per}(\R)$, such that
\begin{equation}\label{eq:PfNonlnrProb_ConvProps4InvDom_inK}
    \u_k\xrightharpoonup{\;\;\;\;\;} \u\;\text{in}\; L^2_\text{per}(\R;\H(\O))\;\;\;\;\;\;\;\text{and}\;\;\;\;\;\;\;\bm{\xi}_k\longrightarrow \bm{\xi}\;\text{in}\; L^2_\text{per}(\R).
\end{equation}

To complete the proof, we must pass to the limit as $k\rw\infty$ in (\ref{eq:LOCALwkForm_coupledSyst_fEQ0Case-finalApprox}). The standard technique uses the convergence properties in (\ref{eq:PfNonlnrProb_ConvProps4InvDom_inK}), but also requires showing the strong convergence
\begin{equation}\label{eq:PfNonlnrProb_StrongCoveAubinLions}
    \u_k\longrightarrow \u\;\text{in}\; L^2_\text{per}(\R;\H(\O_0)).
\end{equation}
For this, one defines the linear functional
\begin{equation*}
        \begin{aligned}
        \left\langle \bm{g}_k, \bm{\Psi} \right\rangle:=\d^2(\bm{F}_0,\bm{\Psi})_{L^2(\O_0)} &-\d\left[ \left((\u_k-\bm{\xi}_k)\cdot\nabla\bm{V}_0,\bm{\Psi}\right)_{L^2(\O_0)}+ \left((\bm{V}_0-\bm{\xi}_0)\cdot\nabla\u_k,\bm{\Psi}\right)_{L^2(\O_0)}\right] \\
        &\;\;\;\;\;\;\;\;\;\;\;\;\;\;\;\;\;\;\;\;\;\;\;\;\;\;\;\;\;- \left((\u_k-\bm{\xi}_k)\cdot\nabla\u_k,\bm{\Psi}\right)_{L^2(\O_0)}-2\nu\left( \u_k,\bm{\Psi}\right)_{\H(\O_0)},
    \end{aligned}
\end{equation*}
for given $\bm{\Psi}\in\H(\O_0)$, and shows that $\bm{g}_k\in L^1_\text{per}(\R;\H(\O_0)')$. It is then easily determined that $\left\langle \frac{\partial \u_k}{\partial t},\bm{\Psi}\right\rangle= \left\langle \bm{g}_k,\bm{\Psi}\right\rangle$, for a.e. $t\in (0,T)$, implying, by the Aubin-Lions Lemma and the series of embeddings\footnote{Here, we use the notation ``$\hookrightarrow\hookrightarrow$" to indicate compact embedding.}
    $$
        \H(\O_0) \hookrightarrow\hookrightarrow \K(\O_0)\hookrightarrow [\H(\O_0)]',
    $$
the validity of (\ref{eq:PfNonlnrProb_StrongCoveAubinLions}), up to a subsequence that may depend on $k_0$. However, with a Cantor diagonalization, one lifts this dependence, so that (\ref{eq:PfNonlnrProb_StrongCoveAubinLions}) holds for all $\O_k$ with $\O_k\supset \text{supp}\,\bm{\f}$. Finally, with (\ref{eq:PfNonlnrProb_ConvProps4InvDom_inK}) and (\ref{eq:PfNonlnrProb_StrongCoveAubinLions}) established, by a routine procedure, one then passes to the limit in (\ref{eq:LOCALwkForm_coupledSyst_fEQ0Case-finalApprox}) as $k\rw\infty$, thereby completing the proof.
\end{proof}

\section{\textbf{Sufficient Conditions for Propulsion}}\label{SCFP}

We are now ready to address the problem of propulsion. Continuing the procedure we started in Section \ref{WkSolNonlnrProb}, we now formally substitute the secondary scaling
\begin{equation}\label{eq:FinalScaling_deltaSqrd_TransOnly}
    \u:=\d^2\w,\;\;\;\;\;\;\;\pi=\d^2\mathsf{p},\;\;\;\;\;\;\;\bm{\xi}:=\d^2\bm{\chi}
\end{equation}
into (\ref{eq:eom_coupledSyst_fEQ0Case-2}), using (\ref{eq:eom_coupledSyst_fEQ0Case-1}) and take the average over $(0,T)$ to obtain
\begin{equation}\label{eq:eom_Splitting_NonlnrEq_Sub2AvgedTransOnly}
\begin{aligned}
    \left.\begin{array}{c}
        \displaystyle
        \d^2\overline{(\w-\bm{\chi})\cdot\nabla\w}+\d \left[ \overline{(\w-\bm{\chi})\cdot\nabla\bm{V}_0}+\overline{(\bm{V}_0-\bm{\xi}_0)\cdot\nabla\w}\right]\\
        \;\;\;\;\;\;\;\;\;\;\;\;\;\;\;\;\;\;\;\;\;\;\;\;\;\;\;\;\;\;\;\;\;\;\;\;\;\;\;\;\;\;\;\;\;\;\;\;\;\;\;\;\;\;\;\;\;\;\;\;\;\;\;\;\;\;\;\;\;\;\;\;\;\;\;\;=\div\,\textbf{T}(\overline{\w},\overline{\mathsf{p}})+\overline{\bm{F}}_0 \\
        \;\;\;\;\;\;\;\;\;\;\;\;\;\;\;\;\;\;\;\;\;\;\;\;\;\;\;\;\;\;\;\;\;\;\;\;\;\;\;\;\;\,\div\,\overline{\w}=0 \\
    \end{array}\right\}&\;\;\;\;\;\;\;{\text{in}\;\O\times\R} \\
    \overline{\w}=\overline{\bm{\chi}}&\;\;\;\;\;\;\;\text{on}\;\partial\O\times\R \\
    \lim_{|\x|\rw\infty}\overline{\w}(\x,t)=\textbf{0}\;&\;\;\;\;\;\;\;\text{in}\;\R. \\
        \displaystyle
        \int_{\partial\O}\textbf{T}(\overline{\w},\overline{\mathsf{p}})\cdot\bm{n}\,\text{d}S=\textbf{0}&\;\;\;\;\;\;\;{\text{in}\;\R}.
\end{aligned}
\end{equation}
Formally taking the limit in (\ref{eq:eom_Splitting_NonlnrEq_Sub2AvgedTransOnly}) as $\delta\rightarrow 0$, we finally arrive at the following steady Stokes problem:
\begin{equation}\label{eq:eom_Splitting_NonlnrEq_Sub2Avged_delta2Zero_TransOnly}
\begin{aligned}
    \left.\begin{array}{c}
        \displaystyle
        \div\,\textbf{T}(\w_0,\mathsf{p}_0)=-\overline{\bm{F}}_0 \\
        \;\;\;\;\;\div\,\w_0=0 \\
    \end{array}\right\}&\;\;\;{\text{in}\;\O}\,, \\
    \w_0=\bm{\chi}_0&\;\;\;\text{on}\;\partial\O\,, \\
    \lim_{|\x|\rw\infty}\w_0(\x)=\textbf{0}\;&\,, \\
        \displaystyle
        \int_{\partial\O}\textbf{T}(\w_0,\mathsf{p}_0)\cdot\bm{n}\,\text{d}S=\textbf{0}&\,.
\end{aligned}
\end{equation}
Dot-multiplying (again, formally) equation (\ref{eq:eom_Splitting_NonlnrEq_Sub2Avged_delta2Zero_TransOnly}$)_1$ by arbitrary $\bm{\psi}\in\H(\Omega)$ and integrating by parts over $\Omega$, as was done to obtain (\ref{eq:wkForm_coupledSyst_fEQ0Case-2}), we are lead to a weak formulation of (\ref{eq:eom_Splitting_NonlnrEq_Sub2Avged_delta2Zero_TransOnly}), made precise by the following definition.

\begin{definition}\label{def:wkSol-2_TransOnly}
    Let $\bm{F}\in L_\text{per}^{\infty}(\mathbb{R})$. Then $(\w_0,\bm{\chi}_0)$ is a {\em weak solution} to the Stokes problem (\ref{eq:eom_Splitting_NonlnrEq_Sub2Avged_delta2Zero_TransOnly}) if
    \begin{enumerate}
        \item $\w_0\in \H(\Omega)$ and $\bm{\chi}_0\in\mathbb{R}^3$ are such that $\w_0=\bm{\chi}_0$ on $\partial\Omega$;
        \item $\w_0$ satisfies
            \begin{equation}\label{eq:wkForm_coupledSyst_fEQ0Case-2_inW_averaged-deltaTakento0_TransOnly}
                2\nu\left( \w_0,\bm{\psi}\right)_{\H}=(\overline{\bm{F}}_0,\bm{\psi})_2,\;\;\;\;\;\text{for every $\bm{\psi}\in \H(\Omega)$}.
            \end{equation}
    \end{enumerate}
\end{definition}

\begin{remark}
    Observe that, because of \eqref{eq:PfExistNonlnrProb_SolonnEmbApplied2NonlnrTerms}$_2$ and  \eqref{GP}, the right-hand side of \eqref{eq:wkForm_coupledSyst_fEQ0Case-2_inW_averaged-deltaTakento0_TransOnly} is well-defined.\label{KTM}
\end{remark}

Now, for each $\d>0$, thanks to Theorem \ref{thm:ExistceWkSol_NonlnrProb}, there exists a corresponding weak solution $(\u_{\delta},\bm{\xi}_{\delta})$ to the nonlinear problem (\ref{eq:eom_coupledSyst_fEQ0Case-2}). We claim that, as $\delta\rightarrow 0$, their scaled time-averaged parts $(\overline{\w}_{\delta},\overline{\bm{\chi}}_{\delta})$ converge to the weak solution $(\w_0,\bm{\chi}_0)$ of (\ref{eq:eom_Splitting_NonlnrEq_Sub2Avged_delta2Zero_TransOnly}), whose existence and uniqueness must be verified first. To this end, we begin by recalling the following result \cite[Section V.4]{galdi2011introduction},   \cite[Sections 5.2-5.4]{happel1983low}.

\begin{lemma}\label{thm:eom_aux_TransOnly}
    Let $s\in(1,\infty)$, $q\in \left(\frac{3}{2},\infty\right)$ and $r\in (3,\infty)$. For each $i=1,2,3$, there exists a unique solution 
$$
    (\bm{h}^{(i)},p^{(i)})\in [D^{2,s}(\Omega)\cap D^{1,q}(\Omega)\cap L^r(\Omega)\cap C^\infty(\Omega)]\times [D^{1,s}(\Omega)\cap L^q(\Omega)\cap C^\infty(\Omega)]
$$ 
to the Stokes problem
    \begin{equation}\label{eq:eom_aux-1}
    \begin{aligned}
    \left.\begin{array}{c}
        \displaystyle
        \div\,\textbf{T}(\bm{h}^{(i)},p^{(i)})=\bm{0} \\
        \;\;\;\;\;\;\;\;\;\;\;\,\div\,\bm{h}^{(i)}=0 \\
    \end{array}\right\}&\;\;\;\;\;\;\;{\text{in}\;\Omega} \\
    \bm{h}^{(i)}=\bm{e}_i&\;\;\;\;\;\;\;\text{on}\;\partial\Omega.
\end{aligned}
\end{equation}
Moreover, for $i,k\in\{1,2,3\}$, we have that the matrix $\mathbb{K}$, defined, component-wise, by
\begin{equation}\label{eq:Def_KMatrix}
    (\mathbb{K})_{ki}:=\bm{e}_k\cdot\int_{\partial\Omega}\left(\textbf{T}(\bm{h}^{(i)},p^{(i)})\cdot\bm{n}\right) \text{d}S
\end{equation}
is both symmetric and invertible.
\end{lemma}

Before we move on to the existence of weak solutions to (\ref{eq:eom_Splitting_NonlnrEq_Sub2Avged_delta2Zero_TransOnly}), we wish to make the following remark about the physical meaning of (\ref{eq:eom_aux-1}). Indeed, this system describes the flow of a viscous liquid around a body with the \textit{prescribed} motion of \textit{pure translation} along basis vector $\bm{e}_i$. In turn, $(\mathbb{K})_{ki}$ is the {\em resistance matrix} and represents the $k^\text{th}$ component of the hydrodynamic force exerted on $\partial\Omega$ due to translational motion along the direction $\bm{e}_i$.

\begin{lemma}\label{thm:StokesProb_ZeroAvgForce}
Problem (\ref{eq:eom_Splitting_NonlnrEq_Sub2Avged_delta2Zero_TransOnly}) admits one and only one weak solution $(\w_0,\bm{\chi}_0)$ such that 
$$
    (\w_0,\bm{\chi}_0)\in [D^{2,s}(\Omega)\cap D^{1,q}(\Omega)\cap L^r(\Omega)]\times \mathbb R^3\,,\ \ s\in(1,\infty), \ q\in(\mbox{$\frac32$},\infty),\ r\in (3,\infty)\,.
$$ 
Furthermore, $\w_0\in C^\infty(\Omega)$, and there exists $\mathsf{p}_0\in C^\infty(\Omega)\cap D^{1,s}(\Omega)\cap L^q(\Omega)$ such that $(\w_0,\mathsf{p}_0,\bm{\chi}_0)$ solves (\ref{eq:eom_Splitting_NonlnrEq_Sub2Avged_delta2Zero_TransOnly}) in the ordinary sense.
\end{lemma}
\begin{proof}
Let
\begin{equation}\label{def:wkSoltoStokesProb-1_TransOnly_fAv_0}
    \bm{\chi}_0=\mathbb{K}^{-1}\cdot\sum_{i=1}^3(\overline{\bm{F}}_0,\bm{h}^{(i)})_2\bm{e}_i\,.
\end{equation}
Since $\bm{h}^{(i)}\in\mathcal H(\Omega)$, $i=1,2,3$, it follows that $\bm{\chi}_0$ is well defined; see Remark \ref{KTM}. With $\bm{\chi}_0=\chi_{0i}\bm{e}_i$, define
\begin{equation}\label{def:wkSoltoStokesProb-2_TransOnly_fAv_0}
    \w_0:=\sum_{i=1}^3\chi_{0i}\bm{h}^{(i)}\;\;\;\;\;\;\;\;\;\;\text{and}\;\;\;\;\;\;\;\;\;\;\mathsf{p}_0:=\sum_{i=1}^3\chi_{0i}p^{(i)}.
\end{equation}
In view of Lemma \ref{thm:eom_aux_TransOnly}, we infer that $(\w_0,\mathsf{p}_0,\bm{\chi}_0)$, possesses all the stated regularity properties. Moreover, multiplying (\ref{eq:eom_aux-1}$)_\text{1-3}$ by $\chi_{0i}$ and summing over $i$, we immediately see that this choice of $(\w_0,\mathsf{p}_0,\bm{\chi}_0)$ satisfies (\ref{eq:eom_Splitting_NonlnrEq_Sub2Avged_delta2Zero_TransOnly}$)_\text{1-3}$. Next, applying $\mathbb{K}$ to the on the left of both sides of (\ref{def:wkSoltoStokesProb-1_TransOnly_fAv_0}), from the definition (\ref{eq:Def_KMatrix}), we also see the validity of (\ref{eq:eom_Splitting_NonlnrEq_Sub2Avged_delta2Zero_TransOnly}$)_4$, completing the proof of the existence.

As for the uniqueness, let $\w_0'$ be another weak solution to (\ref{eq:eom_Splitting_NonlnrEq_Sub2Avged_delta2Zero_TransOnly}) in the sense of (i)-(ii) of Definition \ref{def:wkSol-2_TransOnly}. Then, $\w_0'$ satisfies
$$
    2\nu\left(\w_0',\bm{\psi}\right)_{\H}=\left(\overline{\bm{F}}_0,\bm{\psi}\right)_2,\;\;\;\;\;\text{for every $\bm{\psi}\in \H(\Omega)$}.
$$
    The result then follows by subtracting this from (\ref{eq:wkForm_coupledSyst_fEQ0Case-2_inW_averaged-deltaTakento0_TransOnly}) and taking, in particular, $\bm{\psi}:=\w_0-\w_0'$.
\end{proof}

We are now equipped to prove the convergences claimed above.

\begin{lemma}\label{thm:thrustApprox_fZeroAvgCase_TransOnly}
    Let $\d\in (0,\d_0)$, where $\d_0$ is the positive number given in Theorem \ref{thm:ExistceWkSol_NonlnrProb}, and let $(\u_{\d},\bm{\xi}_{\d})$ be a weak solution to problem (\ref{eq:eom_coupledSyst_fEQ0Case-2}) corresponding to the unique pair $(\bm{V}_0,\bm{\xi}_0)$ given in Theorem \ref{thm:SolLnrProb_TransOnly}. Apply the same scaling in (\ref{eq:FinalScaling_deltaSqrd_TransOnly}) to $\u_{\d}$ and $\bm{\xi}_{\d}$:
\begin{equation}\label{eq:FinalScaling_DeltaSqr_ZeroAvg_TransOnly}
    \u_{\d} :=\d^2\w_{\d},\;\;\;\;\;\;\;\;\;\;\;\bm{\xi}_{\d}:=\d^2\bm{\chi}_{\d}.
\end{equation}
    Then, as $\d\rw 0$, we have
    \begin{equation}\label{eq:PropConvgces_TransOnly}
    \begin{aligned}
        \overline{\w}_{\d} & \xrightharpoonup{\;\;\;\;\,} \w_0\;\;\text{in $\H(\O)$} \\
        \overline{\bm{\chi}}_{\d}&\longrightarrow \bm{\chi}_0\;\;\,\text{in $\R^3$}
    \end{aligned}
    \end{equation}
    where $(\w_0,\bm{\chi}_0)$ is the weak solution to problem (\ref{eq:eom_Splitting_NonlnrEq_Sub2Avged_delta2Zero_TransOnly}) furnished by Lemma \ref{thm:StokesProb_ZeroAvgForce}.
\end{lemma}
\begin{proof}
    By the uniqueness property afforded by Lemma \ref{thm:StokesProb_ZeroAvgForce}, it suffices to show the properties (\ref{eq:PropConvgces_TransOnly}) along a subsequence, say $\{\d_n\}_{n\in\N}$, of positive numbers with $\lim_{n\rw\infty}\d_n=0$. Given such a subsequence, for each $n\in\N$, write
    \begin{equation}\label{eq:scaledQties-2_TransOnly}
        \bm{u}_n=\d_n^2\w_n,\;\;\;\;\;\bm{\xi}_n=\d_n^2\bm{\chi}_n.
    \end{equation}
    Then, from (\ref{eq:UnifEst_WkSol_NonlnrProb_Omega}) and (\ref{eq:scaledQties-2_TransOnly}), we easily obtain the estimate
    \begin{equation}\label{eq:est_on_wn_TransOnly}
        \|\overline{\w}_n\|_{\H}+|\overline{\bm{\chi}}_n|\leq \kappa,
    \end{equation}
    where, from now on, by $\kappa$ we denote a generic positive constant depending, at most, on $T$, $\bm{V}_0$, and $\bm{\xi}_0$. Then, by standard compactness theorems, one finds $\widetilde{\w}\in\H(\O)$ and $\widetilde{\bm{\chi}}\in\R^3$, such that (up to subsequence),
    \begin{equation}\label{eq:propConvergences_TransOnly}
        \overline{\w}_n \xrightharpoonup{\;\;\;\;\,} \widetilde{\w}\;\;\;\text{in $\H(\Omega)$},\;\;\;\;\;\;\;\overline{\bm{\chi}}_n\longrightarrow \widetilde{\bm{\chi}}\;\;\;\text{in $\mathbb{R}^3$}
    \end{equation}
as $n\rw\infty$ and also
    \begin{equation}\label{eq:trace-result-in-prop-thm_TransOnly}
        \widetilde{\w}=\widetilde{\bm{\chi}}\;\;\;\text{on $\partial\Omega$}.
    \end{equation}
Next, in (\ref{eq:wkForm_coupledSyst_fEQ0Case-2}) we take, in particular, $\bm{\f}\in\C(\O)$  and $\bm{u}\equiv\u_n$,  $\bm{\xi}\equiv \bm{\xi}_n$ as in (\ref{eq:scaledQties-2_TransOnly}).  This yields
    \begin{equation}\label{eq:initial_est_wAn_TransOnly}
        2\nu(\overline{\bm{w}}_n,\bm{\f})_{\H(\Omega)}=\d_n^2 A_n+\d_nB_n+(\overline{\bm{F}}_0,\bm{\f})_2,
    \end{equation}
where
\begin{equation}\label{eq:AnBnDef_TransOnly}
\begin{aligned}
    A_n&:=\left(\overline{(\w_n-\bm{\chi}_n)\cdot\nabla\w_n},\bm{\f}\right)_2, \\
    B_n&:=\left[ \left(\overline{(\w_n-\bm{\chi}_n)\cdot\nabla\bm{V}_0},\bm{\f}\right)_2+ \left(\overline{(\bm{V}_0-\bm{\xi}_0)\cdot\nabla\w_n},\bm{\f}\right)_2\right].
\end{aligned}
\end{equation}
Then, comparing (\ref{eq:initial_est_wAn_TransOnly}) with (\ref{eq:wkForm_coupledSyst_fEQ0Case-2_inW_averaged-deltaTakento0_TransOnly}), one immediately sees from (\ref{eq:propConvergences_TransOnly}$)_1$ and (\ref{eq:PropConvgces_TransOnly}$)_1$ that the lemma is proved once we show that $A_n$ and $B_n$ are bounded uniformly in $n$; indeed, then we can pass to the limit in (\ref{eq:initial_est_wAn_TransOnly}) as $n\rightarrow\infty$ and, with (\ref{eq:trace-result-in-prop-thm_TransOnly}), use the uniqueness property of Lemma \ref{thm:StokesProb_ZeroAvgForce}. But this is easily accomplished by first observing, from (\ref{eq:UnifEst_WkSol_NonlnrProb_Omega}) and (\ref{eq:scaledQties-2_TransOnly}), that
\begin{equation}\label{eq:est_for_wn-chin_TransOnly}
    \|\w_n\|_{L^2(\H)}+\|\bm{\chi}_n\|_{L^2(0,T)}\leq \kappa,
\end{equation}
and then employing in (\ref{eq:AnBnDef_TransOnly}) this uniform bound (\ref{eq:est_for_wn-chin_TransOnly}) in combination with Lemma \ref{thm:KornsIdentity_coupledSyst} and H\"older inequality. The proof is thereby complete.
\end{proof}

As a corollary to Lemma \ref{thm:thrustApprox_fZeroAvgCase_TransOnly}, we have now the main result of this section. We recall that
$$
    \bm{F}_0=-(\bm{V}_0-\bm{\xi}_0)\cdot\nabla\bm{V}_0.
$$

\begin{thm}\label{thm:suffCond4prop_ZeroAvgForce}
    Let $(\v,\bm{\g})$ be a weak solution to problem (\ref{eq:eom_coupled_system-TransOnly}) corresponding to the force $\bm{f}\in L_\text{per}^{\infty}(\mathbb{R})$, where $\overline{\bm{f}}=:\delta\, \overline{\bm{F}}=0$. Then, if 
\begin{equation}\label{SP_ZeroAvgForce_TransOnly}
    \bm{G}:=-\sum_{i=1}^3\left(\overline{(\bm{V}_0-\bm{\xi}_0)\cdot\nabla\bm{V}_0},\bm{h}^{(i)}\right)_2\bm{e}_i \neq \textbf{0}\,,
\end{equation} 
we necessarily have $\overline{\bm{\g}}\neq \textbf{0}$; that is, $\mathscr{B}$ experiences propulsion. Precisely, there is $\delta_0>0$ such that
    \begin{equation}\label{*_ZeroAvgForce_TransOnly}
        \overline{\bm{\g}} = \delta^2\mathbb{K}^{-1}\cdot\bm{G}+\bm{\mathsf{R}}(\delta)\,,\ \ \mbox{for all $\delta\in(0,\delta_0)$}\,,
\end{equation}
where
\begin{equation}\label{**_ZeroAvgForce_TransOnly}
    \lim_{\delta\to 0}\frac{1}{\delta^2}\,\bm{\mathsf{R}}(\delta)=0\,.
\end{equation}
\end{thm}
\begin{proof}
    Dot-multiplying both sides of (\ref{eq:eom_aux-1})$_1$ by the difference $\overline{\w}_{\d}-\w_0$ and integrating by parts over $\O$, using (\ref{eq:eom_Splitting_NonlnrEq_Sub2Avged_delta2Zero_TransOnly}$)_3$, (\ref{eq:FinalScaling_DeltaSqr_ZeroAvg_TransOnly}), property (i) of Definition \ref{def:wkSol_NonLrnProbinDelta}, and (\ref{eq:Def_KMatrix}), we get
    \begin{equation*}
        2\nu\sum_{i=1}^3 (\overline{\w}_{\d}-\w_0,\bm{h}^{(i)})_{\H}\bm{e}_i=\mathbb{K}\cdot(\overline{\bm{\chi}}_{\d}-\bm{\chi}_0).
    \end{equation*}
    Taking
    $$
        \bm{\mathsf{R}}(\d):=2\nu\d^2\mathbb{K}^{-1}\cdot\sum_{i=1}^3 (\overline{\w}_{\d}-\w_0,\bm{h}^{(i)})_{\H}\bm{e}_i,
    $$
    this is
    \begin{equation}\label{eq:ExplicitExpress4R_ZeroAvg_TransOnly}
        \bm{\mathsf{R}}(\d)=\d^2(\overline{\bm{\chi}}_{\d}-\bm{\chi}_0),
    \end{equation}
    from which, thanks to (\ref{eq:PropConvgces_TransOnly}$)_2$, property (\ref{**_ZeroAvgForce_TransOnly}) immediately follows. Then (\ref{*_ZeroAvgForce_TransOnly}) is a consequence of applying (\ref{eq:FinalScaling_DeltaSqr_ZeroAvg_TransOnly}$)_2$, (\ref{def:wkSoltoStokesProb-1_TransOnly_fAv_0}), (\ref{eq:decomp_1stappearance}$)_2$, and the fact that $\overline{\bm{\xi}}_0=\textbf{0}$, to (\ref{eq:ExplicitExpress4R_ZeroAvg_TransOnly}).
\end{proof}

\section{\textbf{Numerical Results}}\label{sec:NumericalResults}

In the following section, we numerically investigate the following: if the periodic force $\bm{f}$ acting on the body $\B$ is due to an internal oscillating mass system (see Figure \ref{Schematic}), for which forces $\bm{f}$ and which shapes of $\B$ will $\bm{f}$ result in a non-zero net motion of $\B$? We restrict ourselves to bodies having rotational symmetry around the $z$-axis. We also assume that $\bm{f}$ and $\bm{\gamma}$ are parallel to the $z$-axis. Under these conditions, we seek for solutions possessing the above symmetry, so that the complete numerical setup can be reduced to a two-dimensional setting.

We begin by analyzing the functional value of $\bm{G}$ in (\ref{SP_ZeroAvgForce_TransOnly}) and the linearized periodic system~\eqref{eq:eom_coupledSyst_fEQ0Case-1} for which we must also discretize the stationary Stokes problem~\eqref{eq:eom_aux-1}. Successively, we will consider the fully nonlinear coupled system~\eqref{eq:eom_coupled_system-TransOnly} to  explore whether a non-zero motion can be obtained even when the functional value of $\bm{G}$ is identically zero.

To simplify the numerical analysis, we shall also non-dimensionalize problems (\ref{eq:eom_coupled_system-TransOnly}) and (\ref{eq:eom_coupledSyst_fEQ0Case-1}). However, before doing so, recalling the expressions (\ref{eq:firsAppofForce}) and (\ref{eq:scal_of_f_with_rho}) relating to the force, we first set
$$
    \bm{f}(t):=\mathfrak{m}\ddot{\bm{y}}(t),
$$
where we have ``normalized" the mass $m:=\rho \mathfrak{m}$ and set $\bm{y}:= -y\widehat{\bm{b}}$, for convenience. Then, letting $a$ and $\omega$ denote the characteristic length of $\B$ and the frequency of oscillations, respectively, and setting $\bm{F}:=\mathfrak{m}\ddot{\bm{Y}}$, for $\bm{y}=\d\bm{Y}$, we proceed to scale
\begin{equation}\label{eq:scalings4nonDim}
\begin{aligned}
    \text{$\x$, $\bm{y}$ and $\bm{Y}$} \hspace{0.3cm}&\text{by}\hspace{0.3cm} a, \\
    \text{$t$} \hspace{0.3cm}&\text{by}\hspace{0.3cm} 1/\o, \\
    \text{$\bm{v}$, $\bm{\g}$, $\bm{\xi}_0$, and $\bm{V}_0$} \hspace{0.3cm}&\text{by}\hspace{0.3cm} \o a, \\
    \text{$p$ and $p_0$} \hspace{0.3cm}&\text{by}\hspace{0.3cm} \nu\o, \\
    \text{$M$ and $\mathfrak{m}$} \hspace{0.3cm}&\text{by $a^3$}
\end{aligned}
\end{equation}
to obtain the resulting dimensionless forms
\begin{equation}\label{eq:eom_coupled_system-TransOnly_nonDim}
\begin{aligned}
    \left.\begin{array}{c}
        \displaystyle
        2h^2\left[\frac{\partial\bm{v}}{\partial t}+(\bm{v}-\bm{\g})\cdot\nabla\bm{v}\right]=\div \,\textbf{T}(\bm{v}, p) \\
        \;\;\;\;\;\;\;\;\;\;\;\;\;\;\;\;\;\;\;\;\;\div\,\bm{v}=0 \\
    \end{array}\right\}&\;\;\;\;\;\;\;{\text{in}\;\Omega\times\mathbb{R}} \\ \bm{v}=\bm{\g}&\;\;\;\;\;\;\;\text{on}\;\partial\Omega\times\mathbb{R} \\
    \lim_{|\bm{x}|\rightarrow\infty}\bm{v}(\bm{x},t)=\textbf{0}\;&\;\;\;\;\;\;\;\text{in}\;\mathbb{R}. \\
        \begin{array}{c}
        \displaystyle
        2h^2\dot{\bm{\g}}=2h^2\ddot{\bm{y}}-\int_{\partial\Omega}\textbf{T}(\bm{v}, p)\cdot\bm{n}\,\text{d}S \\
    \end{array}&\;\;\;\;\;\;\;{\text{in}\;\mathbb{R}}
\end{aligned}
\end{equation}

and
\begin{equation}\label{eq:lnzed_syst_non-dim}
\begin{aligned}
    \left.\begin{array}{c}
        \displaystyle
        2h^2\frac{\partial\bm{V}_0}{\partial t}=\div\,\textbf{T}(\bm{V}_0, p_0) \\
        \div\bm{V}_0=0 \\
    \end{array}\right\}&\;\;\;\;\;\;\;{\text{in}\;\O\times\R} \\
    \bm{V}_0=\bm{\xi}_0&\;\;\;\;\;\;\;\text{on}\;\partial\O\times\R \\
    \lim_{|\x|\rw\infty}\bm{V}_0(\x,t)=\textbf{0}\\
    \displaystyle 2h^2\dot{\bm{\xi}}_0=2h^2\ddot{\bm{Y}}-\int_{\partial\O} \textbf{T}(\bm{V}_0, p_0)\cdot\bm{n}\,\text{d}S&\;\;\;\;\;\;\;{\text{in}\;\R},
\end{aligned}
\end{equation}
of (\ref{eq:eom_coupled_system-TransOnly}) and (\ref{eq:eom_coupledSyst_fEQ0Case-1}), respectively, where we have set the dimensionless masses to 1 and defined the dimensionless parameter,
\begin{equation*}\label{eq:def_StkNo}
    h:=a\sqrt{\frac{\o}{2\nu}},
\end{equation*}
representing the, so-called, ``Stokes number". Finally, we similarly non-dimensionlize the auxiliary system (\ref{eq:eom_aux-1}) by scaling~\footnote{Specifically, we are scaling $\bm{h}^{(i)}$ with the characteristic velocity of the system (\ref{eq:eom_aux-1}), which happens to have a value of 1. As a result, the resulting scaled vector is indeed dimensionless.} $\bm{h}^{(i)}$ with 1  and $p^{(i)}$ with $\nu/a$ to obtain the following dimensionless form of the vector $\bm{G}$, defined in (\ref{SP_ZeroAvgForce_TransOnly}):
\begin{equation}\label{eq:ThrustVect}
    \bm{G}=-2h^2\sum_{i=1}^3\left(\int_{\O}\overline{(\bm{V}_0-\bm{\xi}_0)\cdot\nabla\bm{V}_0}\cdot\bm{h}^{{(i)}}\,\text{d}V\right)\bm{e}_i.
\end{equation}

\subsection{Discretization}

The discretization of (\ref{eq:eom_coupled_system-TransOnly_nonDim}), (\ref{eq:lnzed_syst_non-dim}), and (\ref{eq:eom_aux-1}) is fairly standard, but some care is required, as we are looking for periodic solutions of a system that is expected to undergo oscillations of substantial amplitude with a mean velocity that is zero in the linear and close to zero (compared to its amplitude) in the nonlinear case. Even small absolute errors in the amplitude of the oscillatory solution might have a large impact on the resulting average velocity which then gives rise to a substantial relative error. 

To discretize in time, we use the tangent version of the trapezoidal rule for the rigid body system and for the Navier-Stokes equation, as it is sufficiently accurate and introduces little numerical damping. The tangent (midpoint) version is used so that the force acting on the external body is always evaluated inside the intervals $(0,T),\; (T,2T),$ etc. 

We split the period $T$ into $N$ discrete time steps of size $\Delta t=T/N$ and solve the coupled system of Stokes and rigid body motion (\ref{eq:lnzed_syst_non-dim}) in a semi-explicit iteration. We start with
\[
\bm{\g}_n^{(0)} = \bm{\g}_{n-1},\quad
\bm{v}_n^{(0)} = \bm{v}_{n-1}
\]
and iterate for $l=1,2,\dots$ with a relaxation parameter $\omega\in (0,1]$ which is usually set to $\omega=0.8$: 
\[
\begin{aligned}
    \frac{\bm{v}_n^{(l)}-\bm{v}_{n-1}}{\Delta t}&= \frac{1}{2}\div\, \textbf{T}\big(\bm{v}_{n-1}+\bm{v}_n^{(l)},p_n\big), \\
    \div\, \bm{v}_n^{(l)} &=0, \\
    \frac{\widetilde{\bm{\g}}_n^{(l)}-\bm{\g}_{n-1}}{\Delta t}&=\bm{f}(t_{n-\frac{1}{2}}) - \frac{1}{2}\bm{e}_z\int_{\partial\Omega}\textbf{T}\big( \bm{v}_{n-1}+\bm{v}_{n}^{(l-1)},p_n^{(l-1)}\big)\cdot\bm{n}\cdot\bm{e}_z\,\text{d}s,\\
\bm{\g}_n^{(l)} &= \omega\widetilde{\bm{\g}}_n^{(l)}+(1-\omega)\bm{\g}_n^{(l-1)}.
\end{aligned}
\]
By $\widetilde{\bm{\g}}_n^{(l)}$ we denote an intermediate solution that is damped by $\omega$ for better robustness and convergence. Since we only consider setups with cylindrical symmetric and motion only along the axis of symmetry, the Stokes problem is reformulated in cylindrical coordinates reduced to the $(r,z)$-plane. For numerical simulation, we have to restrict the domain $\O=\mathbb{R}^3\setminus\overline{\B}$ by introducing, for each $R>R_*$, the numerical domain $D_{R}$ as
\[
    D_{R} = \left\{ (r,z)\in\mathbb{R}^2\,:\, r\in (0,R),\; z\in (-R,R)\} \setminus \overline{\B}\right\},
\]
where the rigid body $\B$ is centered at $(0,0)$. We take $R=16$. This choice has been shown to have minimal impact on the artificial domain restriction. 

On the symmetry boundary at $r=0$ we impose the Dirichlet condition $v_r=0$ for the radial component. On the outer boundary at $r=R$ a no-slip condition $\bm{v}=\textbf{0}$ is set. On the surface of the body $\partial\B$ the fluid's velocity matches that of the solid (i.e. $\bm{v}=\bm{\g}$). 

For the linear test cases aiming at the functional values, the do-nothing outflow condition $\partial_n \bm{v}-p\bm{n}=\textbf{0}$ holds, see \cite{HeywoodRannacherTurek1992} and is imposed at $z=R$ and $z=-R$. This condition however is not compatible with the transformed body-centered reference system of the nonlinear problem as $\bm{v}$ is not zero on the numerical boundary. 

The surface force $\textbf{T}(\bm{v},p)\cdot\bm{n}$ on the body is evaluated as a volume integral
\[
    \int_{\partial\Omega}\textbf{T}(\bm{v},p)\cdot\bm{n}\cdot\bm{e}_z\,\text{d}s = \int_\Omega \frac{d\bm{v}}{dt}\cdot \bm{\overline{z}}+\textbf{T}(\bm{v},p):\nabla\bm{\overline{z}}\,\text{d}x,
\] 
where $\bm{\overline{z}}\in H^1(\Omega)$ with $\bm{\overline{z}} = \bm{e}_z$ on $\partial\B$ and $\bm{\overline{z}}=\textbf{0}$ on the other boundaries. In comparison with the evaluation of the surface integral, this approach gives higher order convergence, see~\cite{BraackRichter2006}.

In space we discretize with quadratic finite elements for velocity and pressure using the local projection stabilization technique to cope with the missing inf-sup stability, see~\cite{BeckerBraack2001}. The nonlinear algebraic problems are solved by Newton's method and the linear systems by GMRES, using a parallel geometric multigrid iteration as preconditioner. Details on the discretization are given in~\cite[Chapter 4]{Richter2017} and the nonlinear and linear solver is detailed in~\cite{FailerRichter2019}. The numerical mesh is adjusted to the body and the mesh size $\Delta x$ is adjusted to resolve the flow close to the body. On a mesh with average mesh size $\Delta x=0.25$, this size varies from $\Delta x\approx 0.0073$ at the body to $\Delta x\approx 0.5$ on the outer boundaries. 

In dealing with the linear Stokes problem, we know that the resulting periodic solution of (\ref{eq:lnzed_syst_non-dim}) will result in zero average velocities $\bm{\overline{v}}=\textbf{0}$ and $\bm{\overline{\gamma}}=\textbf{0}$. Hence, to accelerate convergence towards the periodic solution we use the following simple algorithm that aims at projecting the average to the expected one (see \cite{Richter2020} for a detailed description of the algorithm and an extension to the time-periodic Navier-Stokes problem). Starting with $\bm{v}_0:=\textbf{0}$ and $\bm{\gamma}_0:=\textbf{0}$ we iterate for $l=1,2,\dots$,
\begin{enumerate}
	\item Approximate the coupled Stokes problem on $I=[0,T]$ using the initial values 
\[
\bm{v}^{(l)}(0) = \bm{v}^{(l-1)}_0,\quad 
\bm{\gamma}^{(l)}(0) = \bm{\gamma}^{(l-1)}_0;
\]
	\item If periodicity is reached up a given tolerance $\epsilon_P>0$, namely
\[
\|\bm{v}^{(l)}(T)-\bm{v}^{(l)}(0)\|+
\|\bm{\gamma}^{(l)}(T)-\bm{\gamma}^{(l)}(0)\| < \epsilon_P,
\]
stop and accept $\bm{v}^{(l)},\bm{\gamma}^{(l)}$;
\item Correct the average and predict new initial value:
\[
\bm{v}^{(l)}_0:=\bm{v}^{(l)}(T) - \frac{1}{T}\int_0^T \bm{v}^{(l)}(t)\,\text{d}t,\quad
\bm{\gamma}^{(l)}_0:=\bm{\gamma}^{(l)}(T) - \frac{1}{T}\int_0^T \bm{\gamma}^{(l)}(t)\,\text{d}t. 
\]
\end{enumerate}

This iteration quickly converges to the periodic solution and the periodic tolerance is usually reached within a few cycles. For the nonlinear coupled problem, we apply this algorithm in the very first step to correct the initial values. 

\subsection{Definition of the Test Cases}

\begin{figure}[t]
\begin{center}
\includegraphics[width=\textwidth]{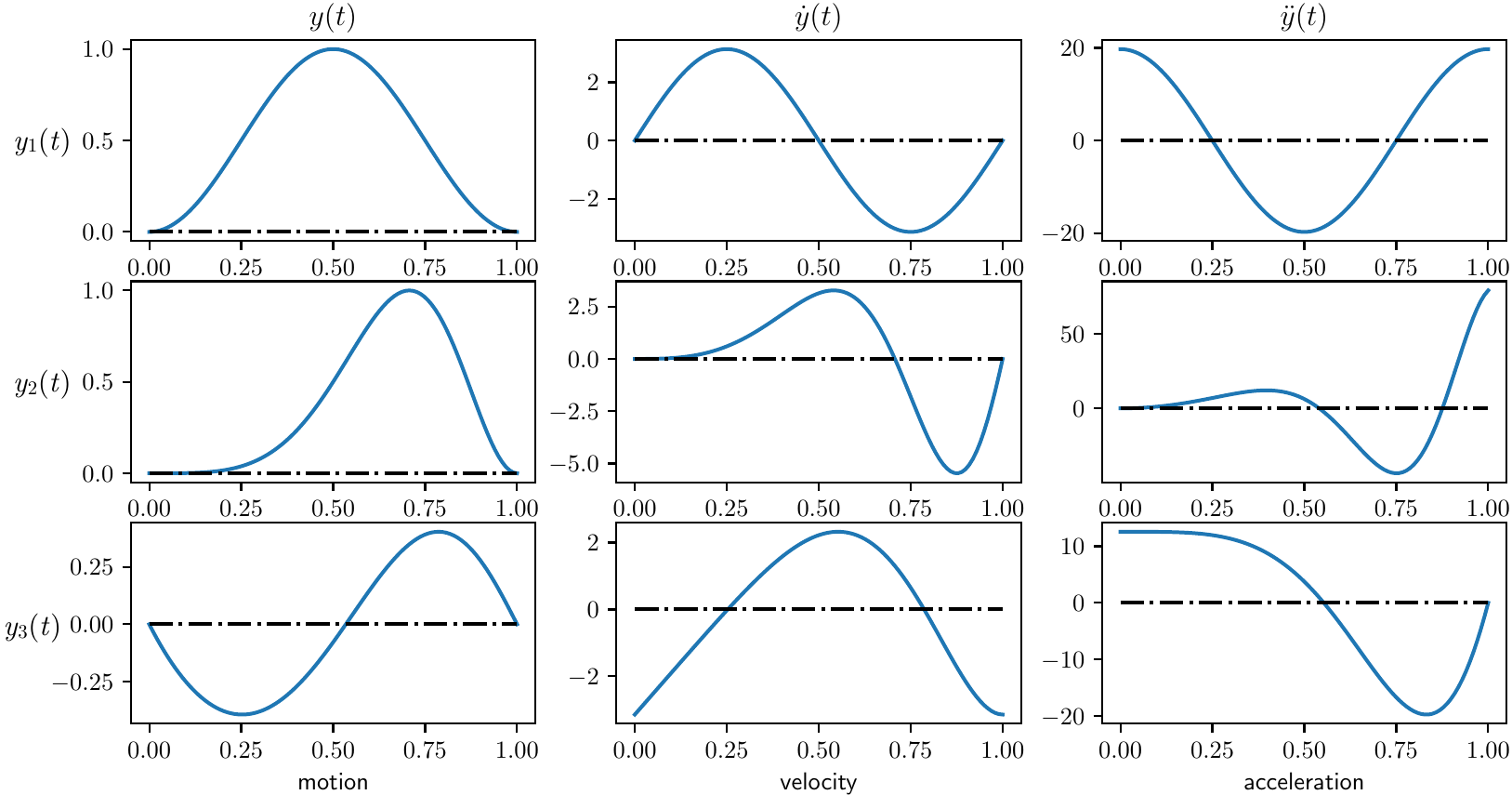}	
\end{center}
\caption{Three forces used in the simulation. From left to right: internal motion of the rigid body $y(t)$, relative velocity of the rigid body $\dot{y}(t)$ and its acceleration $\ddot{y}(t)$. From top to bottom: symmetric smooth force $\ddot{y}(t)$, non-symmetric forces $\ddot{y}_2(t)$ and $\ddot{y}_3(t)$.}
\label{fig:forces}
\end{figure}

Recall that we construct the force acting on $\B$ via the internal motion of a mass. As such, we consider three test forces which we identify with the following (dimensionless) accelerations (see (\ref{eq:eom_coupled_system-TransOnly_nonDim}$)_5$ and (\ref{eq:lnzed_syst_non-dim}$)_5$):
\begin{equation}\label{forces}
\begin{aligned}
    \ddot{y}_1(t)=\ddot{Y}_1(t) &= \frac{\text{d}^2}{\text{d}t^2}\left[\sin^2(\pi t)\right] \\
    \ddot{y}_2(t)=\ddot{Y}_2(t) &= \frac{\text{d}^2}{\text{d}t^2}\left[\sin^2(\pi t^2)\right] \\
    \ddot{y}_3(t)=\ddot{Y}_3(t) &= \frac{\text{d}^2}{\text{d}t^2}\left[\sin(\pi t^2)-t\left(1-t\right)\pi\right].
\end{aligned}
\end{equation}
As such, $\ddot{y}_1$ represents a ``fully symmetric" force with high regularity, whereas $\ddot{y}_2$ and $\ddot{y}_3$ are ``non-symmetric" (see Figure \ref{fig:forces}). Each force is periodically extended from $[0,T)$ to $\mathbb{R}$. As such, $\ddot{y}_1$ is $C^\infty(\mathbb{R})$ and $\ddot{y}_2$ and $\ddot{y}_3$ have anti-derivatives $\dot{y}_2$ and $\dot{y}_3$ that are continuous on $[0,T)$ (see also Figure \ref{fig:forces}).

\begin{figure}[t]
\begin{center}
\includegraphics[width=0.6\textwidth]{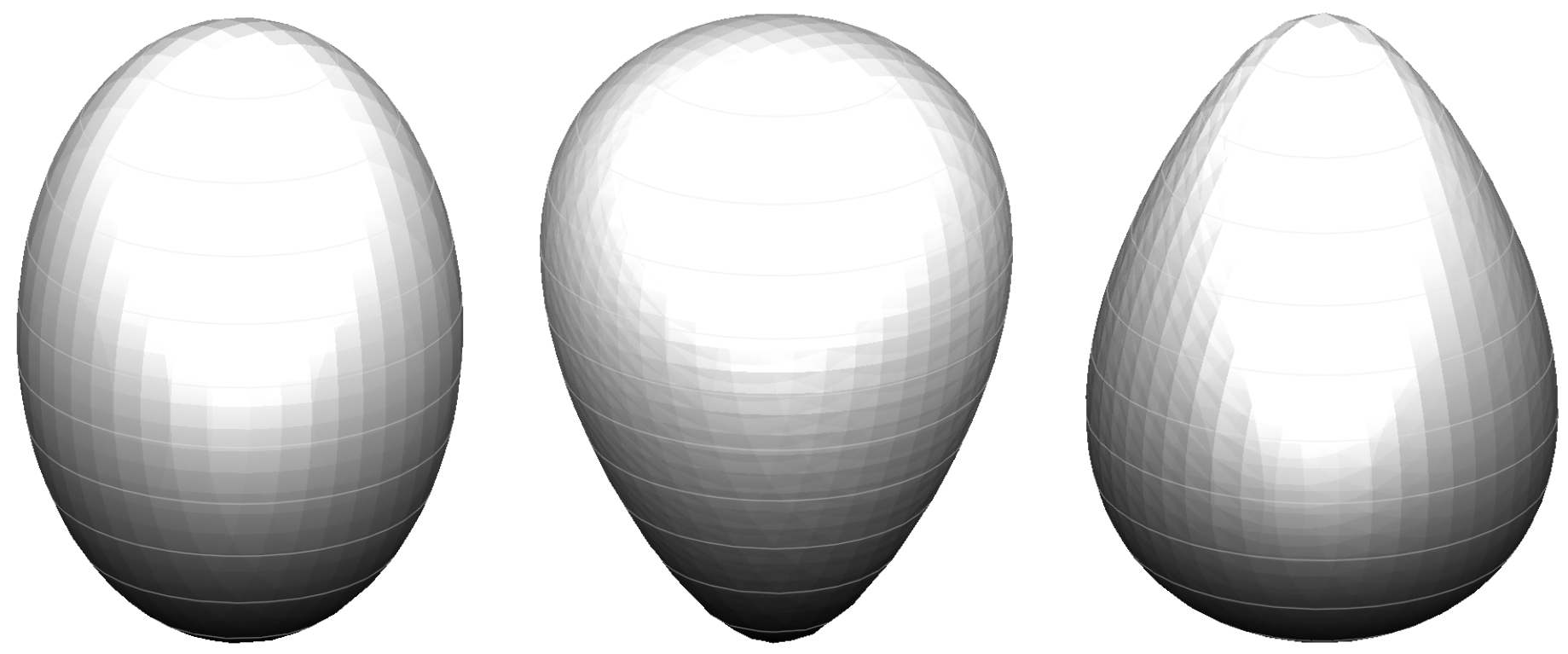}	
\end{center}
\caption{Ellipsoid with symmetry in the direction of flow (left), drop-like shape (center), and flipped drop (right).}
\label{fig:shapes}
\end{figure}

Although all all forces $\ddot{y}_1,\ddot{y}_2$ and $\ddot{y}_3$ have average zero, numerical quadrature on the time grid $0=t_0<\dots<t_N=T$ with the mid point rule introduces an error
\[
    0 \,=\, \int_0^T \ddot{y}(t)\,\text{d}t\,\approx\, \sum_{i=1}^N\Delta t\cdot  \ddot{y}\left(\frac{t_{i-1}+t_i}{2}\right)\,\neq\, 0.
\]
For numerical simulation we therefore correct the force such that its discrete average is exactly zero independent of the step size $\Delta t>0$. Instead of $\ddot{y}(t)$ we use
\[
    \ddot{y}(t)-\frac{1}{T}\sum_{i=1}^N\Delta t\cdot  \ddot{y}\left(\frac{t_{i-1}+t_i}{2}\right)
\]
We study three different bodies, defined as follows:
\[
    \begin{aligned}
    \B_e &= \{ \bm{x}\in \mathbb{R}^3,\; 1.5(x_1^2+x_2^2)+0.7 x_3^2 <1 \},\\
    \B_d &= \left\{ \bm{x}\in \mathbb{R}^3,\; 1.5\frac{x_1^2+x_2^2}{(1+0.3 x_3)^2}+0.7 x_3^2 <1 \right\},\\
    \B_{fd} &= \left\{ \bm{x}\in \mathbb{R}^3,\; 1.5\frac{x_1^2+x_2^2}{(1-0.3 x_3)^2}+0.7 x_3^2 <1 \right\}.
\end{aligned}
\]
Here, $\B_e$ is an ellipsoid, $\B_d$ is drop-like shape, and $\B_{fd}$ is the body $\B_d$ in the opposite orientation with respect to the axis of symmetry (see Figure \ref{fig:shapes}).

\subsection{Computation of the Thrust Vector $\textit{\textbf{G}}$}

We study the linearized problem (\ref{eq:lnzed_syst_non-dim}). In Table~\ref{tab1}, we indicate the functional values of $\bm{G}$, as given in (\ref{eq:ThrustVect}), for the three types of forces and three shapes. The functional value is non-zero only in the case of the non-symmetric bodies $\B_d$ and $\B_{fd}$. The combination of the symmetric force together with a non-symmetric body does generate a non-zero functional value. The value of $\bm{G}$ changes its sign for the mirrored body $\B_{fd}$. It is noteworthy that the three forces $\ddot{y}_1$, $\ddot{y}_2$ and $\ddot{y}_3$  do not give substantially different values.  

\begin{table}[H]
\begin{center}
\begin{tabular}{r|rlll}
\toprule
$h=8$&& $\ddot{y}_1$ & $\ddot{y}_2$ & $\ddot{y}_3$ \\
\midrule
ellipsoid $\B_e$&$\bm{G}=\!\!\!\!$ & \phantom{-}0 & \phantom{-}0 & \phantom{-}0     \\
drop    $\B_d$&$\bm{G}=\!\!\!\!$ & $-0.01124$ & $-0.03943$  & $-0.008484 $\\
flipped drop    $\B_{fd}$&$\bm{G}=\!\!\!\!$ & \phantom{-}0.01124 &  \phantom{-}0.03943&  \phantom{-}0.008484\\
\bottomrule
\end{tabular}

\end{center}
\caption{Computed values of $\bm{G}$ for Stokes number $h=8$ depending on the shape of the domain (elliptic $\B_e$ or drop-like $\B_d$ and $\B_{fd}$) and the type of force (symmetric $\ddot{y}_1$, non-symmetric $\ddot{y}_2$ and $\ddot{y}_3$). As expected, the values of $\bm G$ for $\mathscr B_{d}$ and $\mathscr B_{fd}$ are opposite to each other.}
\label{tab1}
\end{table}

Next, we study the dependence on the Stokes number $h$ of the average body velocity $\overline{\bm{\gamma}}_0:=\mathbb K^{-1}\cdot\bm{G}$ determined in Theorem \ref{thm:suffCond4prop_ZeroAvgForce}. We limit the study to the case of the flipped drop-shaped body $\B_{fd}$ forced by $\ddot{y}_2$. Numerical values and their corresponding plot are reported at the top and left of Table~\ref{tab:frequency}. These results present interesting features. In the first place, $\overline{\bm{\gamma}}_0$ has a quadratic behavior for sufficiently large values of $h$, which implies that, for such $h$, the speed of the body is proportional to the frequency of oscillations. Moreover, for small $h$, there is a ``critical" value of $h$, between $h=3$ and $h=4$, around which $\overline{\bm{\gamma}}_0$ gets smaller before growing quadratically. For comparison we show the resulting average velocity $\overline{\bm{\gamma}}$ of the full nonlinear problem. These values are larger than the corresponding results for $\overline{\bm{\gamma}}_0$ but they match its trajectory for larger values of $h$.\looseness=-1

\begin{table}[H]
\begin{center}
\begin{tabular}{l|lllllllll}
\toprule

$h$ & $1$ & $2$ & $4$ & $8$ & $16$ & $32$ & $64$ & $128$ & 256 \\
\midrule 
$\overline{\bm{\gamma}}_0$&0.00050&
0.0010&
0.00076&
0.039&
0.83&
7.82&
42.62&
177.80&        
713.49\\
\bottomrule 
\end{tabular}

\centerline{\includegraphics[width=0.6\textwidth]{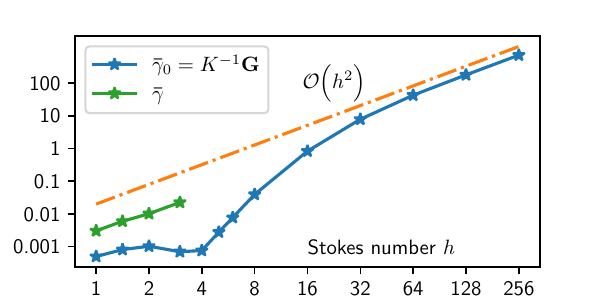}}

\caption{Top: Dependency of $\overline{\bm{\gamma}}_0\equiv {\mathbb K}^{-1}\cdot\bm{G}$ on the Stokes number for the flipped drop $B_{fd}$ with $\ddot{y}_2$. Bottom: Visualization of the dependency compared to the average velocity $\overline{\bm{\gamma}}$ of the corresponding full nonlinear Navier-Stokes problem.}\label{tab:frequency}
\end{center}
\end{table}

\vspace{-24pt}

\subsection{Numerical Integration of the Nonlinear Coupled Problem}

The second-order analysis shows that the non-symmetry of the body yields a non-zero functional value of thrust vector (\ref{eq:ThrustVect}) implying a non-zero net motion of the body. For the (round) ellipsoid however, the functional output is zero, regardless of the choice of force. In order to gain insight into the role of symmetry breaking of body and force, we provide numerical results for the full nonlinear coupled problem. 

All tests are run at Stokes number $h=3$.
We consider the same nine combinations tested previously for $\bm{G}$ in Table \ref{tab1} and report the values of the average velocity of the body $\bm{\overline{\gamma}}$ resulting, this time, directly from the numerical integration of the coupled nonlinear problem \eqref{eq:eom_coupled_system-TransOnly}; the results are given in Table~\ref{tab1-nonlin}. They show, in particular, that, in contrast to the second-order prediction, a non-zero net motion occurs in the case of the ellipsoid whenever the force is not symmetric. In the case of the symmetric force $\ddot{y}_1$, the resulting average velocity for the drop and the flipped drop have opposite sign and the same magnitude up to the numerical discretization error. This is expected due to the symmetry of the problem setup. Studying the two non symmetric forces $\ddot{y}_2$ and $\ddot{y}_3$ reveals the distinct influence of the symmetry of force and shape on the average velocity. In both cases, the mean of the velocities of drop and flipped drop is close to the velocity of the ellipsoid which can be considered as the shape resulting  from overlaying $\B_d$ and $\B_{fd}$.\looseness=-1

\begin{table}[h]
\begin{center}
\begin{tabular}{r|rlll}
\toprule
$h=3$&& $\ddot{y}_1$ & $\ddot{y}_2$ & $\ddot{y}_3$ \\
\midrule
ellipsoid $\B_e$&$\bm{\overline{\gamma}}=\!\!\!\!$ & \phantom{-}0 & 0.0256 & 0.00283 \\
drop   $\B_d$&$\bm{\overline{\gamma}}=\!\!\!\!$ & \phantom{-}0.00677 & 0.0391 & 0.00591\\
flipped drop    $\B_{fd}$&$\bm{\overline{\gamma}}=\!\!\!\!$ &-0.00672 & 0.0225 & 0.00119 \\
\bottomrule
\end{tabular}
\end{center}
\caption{Average body velocity $\bm{\overline{\gamma}}$ computed from the full nonlinear problem for Stokes number $h=3$ and the three forces and shapes. It is interesting to observe that the absolute values of $\overline{\bm{\gamma}}$ for $\mathscr B_{d}$ are  larger than those for $\mathscr B_{fd}$,  indicating, as expected, that the viscous drag on $\mathscr B_{fd}$ is larger than the one  on $\mathscr B_{d}$. }
\label{tab1-nonlin}
\end{table}

In Table~\ref{tab45}, we report the resulting average body velocity $\bm{\overline{\gamma}}$ depending on the Stokes number $h$. We find approximately quadratic growth of the average velocity with increasing Stokes number. This is the same behavior observed for the thrust vector $\bm{G}$. 

\begin{table}[h]
\begin{center}
	\begin{minipage}{0.49\textwidth}
		\includegraphics[width=\textwidth]{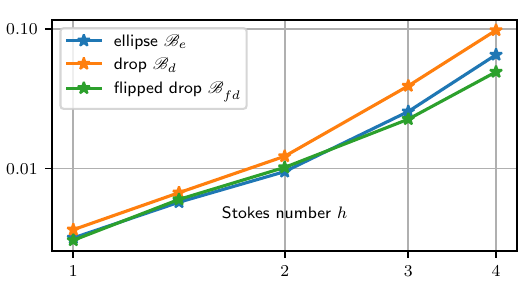}
	\end{minipage}\hspace{0.02\textwidth}
    \begin{minipage}{0.45\textwidth}
		
		\begin{tabular}{l|rrrr}
		\toprule 
		 	$h$ & $1$ & $2$ & $3$ &$4$\\
			\midrule
			$\B_e$   & 0.00317 & 0.0095 & 0.0256 &--\\
			$\B_d$   & 0.00363 & 0.0122 & 0.0391&0.0979  \\
			$\B_{fd}$& 0.00305 & 0.0102 & 0.0225&--\\
			\bottomrule
		\end{tabular}

	\end{minipage}
\end{center}
\caption{Dependency on the Stokes number $h$ of the average body velocity $\bm{\overline{\gamma}}$ computed from the full nonlinear problem, for the three different shapes. All problems are driven with $\ddot{y}_2$.}
\label{tab45}
\end{table}

\begin{figure}[H]
  \centering
  \includegraphics[width=0.65\linewidth]{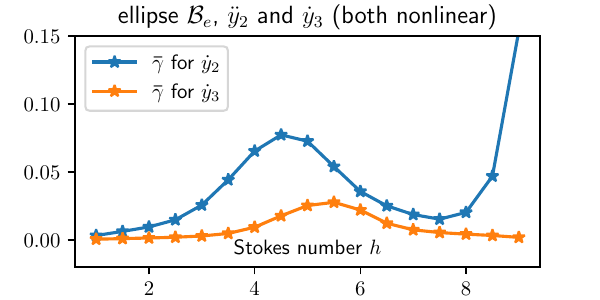}
  \caption{Variation with Stokes number of the average velocity, $\overline{\gamma}$, for the ellipsoid, computed from the full nonlinear problem for forces $\ddot{y}_2$ and $\ddot{y}_3$.}
\label{fig:ns-ellipse-23}
\end{figure}

Figure \ref{fig:ns-ellipse-23} refers to the case where $\mathscr B$ is an ellipsoid. It shows the behavior of the average velocity computed from the full nonlinear problem for different values of the Stokes numbers and for forces $\ddot{y_2}$ and $\ddot{y_3}$. We recall that, as shown in Table \ref{tab1}, the $O(\delta^2)$ theory is not able to capture propulsion in these cases.

\begin{figure}[H]
  \centering
  \includegraphics[width=1\linewidth]{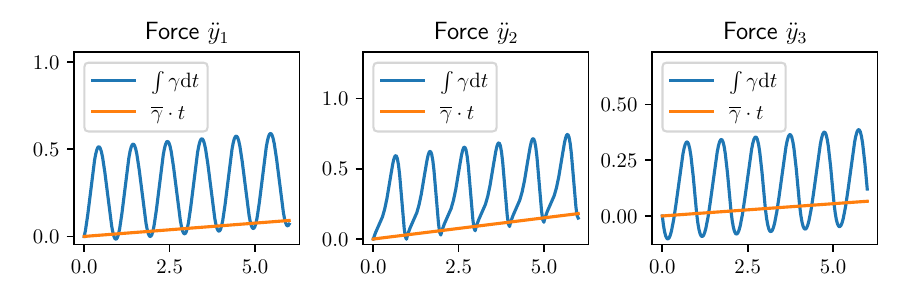}
  \caption{Position of the center of mass of the flipped drop $\mathscr B_{fd}$ ($\int_0^t\gamma(s)\,{\rm d}s$) and net distance covered $(\overline{\gamma}\,t$) vs. time, computed from the full nonlinear problem. }
  \label{traj}
\end{figure}

\noindent
Finally, in Figure \ref{traj} we report, in the case of the flipped drop and for different forces, the variation with time of the position of the center of mass along with its net distance traveled, computed from the full nonlinear problem.

\section*{Acknowledgments}
The work of J. Edelmann and T. Richter has been supported by the German Research Foundation, Grant  537063406, and that of G.P. Galdi and M.M. Karakouzian by US National Science Foundation, Grant DMS 2307811.

%\bibliographystyle{ieeetr} %alpha, apalike, ieeetr
%\bibliography{References.bib}

\end{document}